\documentstyle[12pt]{article}     
\begin{document} 
\title{Stochastic processes on non-Archimedean 
spaces. I. Stochastic processes on Banach spaces.}  
\author{S.V. Ludkovsky}  
\date{05 March 2001
\thanks{Mathematics subject classification
(1991 Revision) 28C20 and 46S10.} }
\maketitle
\par  address: Laboratoire de Math\`ematiques Pures, 
\par Complexe Scientifique des C\`ezeaux,
\par 63177 AUBI\`ERE Cedex, France.
\par permanent address: Theoretical Department,
Institute of General Physics,
Str. Vavilov 38, Moscow, 119991 GSP-1, Russia.
\begin{abstract} 
Non-Archimedean analogs of Markov quasimeasures 
and stochastic processes are investigated. They are used for the 
development of stochastic antiderivations.
The non-Archimedean analog of the It$\hat o$ formula is proved.
\end{abstract} 
\section{Introduction.}
Stochastic differential equations on real Banach spaces and manifolds
are widely used for solutions of mathematical and physical
problems and for construction and investigation of measures
on them \cite{dalf,gihsko,ikwat,mckean,malli,oksen}.
In particular stochastic equations can be used for 
the constructions of quasi-invariant measures on topological groups.
In the cases of real Banach-Lie groups, some simplest
cases of diffeomorphisms groups and free loop spaces
of real manifolds such stochastic equations and measures
were investigated in \cite{aiel,dalschn,eljm,malli}.
Stochastic processes on geometric loop groups and diffeomorphism groups
of wide classes of real and complex manifolds were 
investigated in \cite{lustrcm}.
On the other hand, non-Archimedean functional analysis
develops fastly in recent years and also its applications
in mathematical physics \cite{ami,esc,sch1,sch2,roo,vla3,khrif,jang}.
Wide classes of quasi-invariant measures including 
analogous to Gaussian type
on non-Archimedean Banach spaces, loops and diffeomorphisms 
groups were investigated in \cite{lud,lu6,luumnls,lutmf99,lubp2}.
Quasi-invariant measures on topological groups and their
configuration spaces can be used 
for the investigations of their unitary representations
(see \cite{luseamb,lubp99,lutmf99,lubp2} and references therein).
\par In view of this developments non-Archimedean analogs
of stochastic equations and diffusion processes 
need to be investigated. Some steps in this direction were 
made in \cite{bikvol,evans}. 
There are different variants for such activity, for example,
$p$-adic parameters analogous to time, but spaces
of complex-valued functions. At the same time measures may be 
real, complex or with values in a non-Archimedean field.
\par In the classical stochastic analysis indefinite integrals
are widely used, but in the non-Archimedean case the field of $p$-adic
numbers $\bf Q_p$ has not linear order structure apart from $\bf R$. 
For elements
$f$ in the space of $m-1$ times continuously differentiable functions
$C^{m-1}$ there are antiderivation operators
$P_m: C^{m-1}\to C^m$ such that $(P_mf)'=f$ \cite{sch1,sch2}. Therefore,
in this paper these indefinite integrals are used,
but they are transformed for more complicated needed cases.
In the classical case for the investigations of stochastic processes
nuclear, Hilbert-Schmidt and of the class $L_q$ operators
are used \cite{dalf,pietsch}. In the non-Archimedean case
the operator theory differs from that of classical and the corresponding
definitions and propositions was necessary to give anew in this article.
\par This work treats the case which was not considered by 
another authors and that is suitable and helpful for 
the investigation of stochastic processes and quasi-invariant measures
on non-Archimedean topological groups. These investigations are not restricted
by the rigid geometry class \cite{freput}, since it is rather narrow.
Wider classes of functions and manifolds are considered. This is posssible
with the use of
Schikhof's works on classes of functions $C^n$ in the sence of
difference quotients, which he investigated few years later the
published formalism of the rigid geometry.
Here are considered spaces of functions with values in Banach spaces
over non-Archimedean local fields, in particular, with values in
the field $\bf Q_p$ of $p$-adic numbers.
For this
non-Archimedean analogs of stochastic processes are considered
on spaces of functions with values in the non-Archimedean 
infinite field with a non-trivial valuation such that a parameter 
analogous to the time is $p$-adic (see \S \S 4.1, 4.2). 
Their existence is proved in Theorem 4.3.
Specific antiderivation operators generalizing Schikhof antiderivation 
operators on spaces of functions $C^n$ are investigated (see \S 2).
Their continuity and differentiability properties are given in Lemmas
2.2, 2.3 and Theorem 2.14. Also operators analogous to nuclear operators
are studied (see Definition 2.10 and Propositions 2.11, 2.12).
In \S 3 non-Archimedean analogs of Markov quasimeasures are defined
and Propositions 3.3.1 and 3.3.2 about their boundedness
and unboundedness are proved. The non-Archimedean stochastic integral 
is defined in \S 4.4. Its continuity as the operator on the corresponding 
spaces of functions is proved in Proposition 4.5.
In Theorems 4.6, 4.8 and Corollary  4.7 analogs of the It\^o formula  
are proved.
Spaces of analytic functions lead to simpler expressions of 
the It\^o formula analog, but the space of analytic functions 
is very narrow and though it is helpful in non-Archimedean 
mathematical physics it is insufficient for solutions of all
mathematical and physical problems. 
For example, in many cases of topological groups for non-Archimedean
manifolds spaces of analytic functions are insufficient.
On the other hand, for spaces $C^n$ rather simple formulas are found.
This work was started five years ago, but because of lack of free 
time it was finished only recently.
All results of this paper are obtained for the first time.
\section{Specific antiderivations of operators.}
{\bf 2.1.} Let $X:=c_0(\alpha ,{\bf K_p})$ 
be a Banach space over a local field (see \cite{wei})
$\bf K_p$ such that ${\bf K_p}\supset \bf Q_p$, 
$\{ e_j: j\in \alpha \} $ denotes the 
standard orthonormal base in $c_0(\alpha ,{\bf K_p})$
where $\alpha $ is an ordinal \cite{eng}, $e_j=(0,...,0,1,0,...)$ 
with the unit on the $j$-th place, $j\in \alpha $ \cite{roo}.
The space $c_0(\alpha ,{\bf K_p})$ consists of vectors
$x=(x_j: x_j\in {\bf K_p}, j\in \alpha )$ such that
for each $\epsilon >0$ a set $\{ j: j\in \alpha ; |x_j|>\epsilon \} $
is finite. The norm in it is the following: $\| x \| :=\sup_j|x_j|$.
It is convenient to supply the set $\alpha $ with 
the ordinal structure due to the Kuratwoski-Zorn lemma. Let 
$F$ be a continuous function on $B_r\times C^0(B_r,X)^{\otimes k}$ 
with values in $C^0(B_r,X)$:
$$(1)\mbox{ }F\in C^0(B_r\times C^0(B_r,X)^{\otimes k}, C^0(B_r,X)),$$
where $Z^{\otimes k}=Z\otimes ... \otimes Z$ is the product of $k$ copies 
of a normed space $Z$ and $Z^{\otimes k}$ is supplied with the box 
(maximum norm) topology \cite{eng,nari},
$B_r:=B({\bf K_p},t_0,r)$ is a ball in $\bf K_p$
containing $t_0$ and of radius $r$, Banach spaces $C^t(M,X)$
of mappings $f: M\to X$ from a $C^{\infty }$-manifold $M$ with clopen charts
modelled on a Banach space $Y$ over $\bf K_p$
into $X$ of class of smoothness $C^t$ with $0\le t < \infty $
are the same as in \cite{luseamb,lutmf99,ludgla,lubp2} with the supremum-norm,
when $M$ is closed and bounded in the corresponding Banach space.
Such mappings can be written in the following form:
$$(2)\mbox{ }F(v,\xi _1,...\xi _l)=\sum_{j\in \alpha }
F^j(v,\xi _1,...,\xi _k)e_j,$$ 
where $F^j\in C^0(B_r\times C^0(B_r,X)^{\otimes k},{\bf K_p})$
for each $j\in \alpha $.
In particular let 
$$(3)\mbox{ }F(v;\xi _1,...,\xi _k)=G(v;\xi _1,...,\xi _l).
(A_{l+1}(v)\xi _{l+1},...,A_k(v)\xi _k),$$
where  $\mbox{ }L(X,Y)$ denotes a Banach space of continuos linear operators 
$A: X\to Y$ supplied with the operator norm 
$\| A \| :=\sup_{0\ne x\in X} \| Ax \| _Y/ \| x \| _X$
and $L(X):=L(X,X)$,
$A_i(v)$ are continuous linear operators
for each $v\in B_r$ such that $A_i\in C^0(B_r,L(X))$, 
$G(v,\xi _1,...,\xi _l)\in L_{k-l}(X^{\otimes (k-l)};X)$
for each fixed $v\in B_r$ and $\xi _1,...\xi _l \in C^0(B_r,X)$,
that is, $F$ is a $(k-l)$-linear operator by $\xi _{l+1},....,\xi _k$,
where $G=G(v,\xi _1,...,\xi _l)$ is the short notation of
$G(v,\xi _1(v),...,\xi _l(v)),$ $L_k(X_1,...,X_k;Y)$ denotes
the Banach (normed) space of $k$-linear continuous operators from 
$X_1\otimes ... \otimes X_k$ into $Y$ for Banach (normed) spaces
$X_1,...,X_k, Y$ over $\bf K$ and $L_k(X^{\otimes k};Y):=L_k(X_1,...,X_k;Y)$
for the particular case $X_1=...=X_k=X$. When $l=0$ put $G=G(v)$.
There exists the following antiderivation of operators given by equation $(3)$:
$$(4)\mbox{ }{\hat P}_{(\xi _{l+1},...,\xi _k)}[G(s;\xi _1,...,
\xi _l)\circ (A_{l+1}\otimes ...\otimes A_k)](v):=$$
$$\sum_{n=0}^{\infty }G(v_n;\xi _1,...,\xi _l).(A_{l+1}(v_n)[\xi _{l+1}(
v_{n+1})-\xi _{l+1}(v_n)],...,A_k(v_n)[\xi _k(v_{n+1})-\xi _k(v_n)]),$$
where $v_n=\sigma _n(t),$, $\{ \sigma _n: n=0,1,2,.. \} $
is an approximation of the identity in $B_r$. 
By its definition the approximation 
of the identity satisfies the following conditions:
\par $(i)$ $\sigma _0(t)=t_0,$
\par $(ii)\mbox{ }\sigma _m\circ \sigma _n=\sigma _n\circ \sigma _m
\mbox{ for each }m\ge n$
and there exists $0<\rho <1$ such that from
\par $(iii)\mbox{ }|x-y|<\rho ^n\mbox{ it follows } \sigma _n(x)=\sigma _n(y)
\mbox{ and }$ 
\par $(iv)\mbox{ }|\sigma _n(x)-x|<\rho ^n$
(see \S 62 and \S 79 in \cite{sch1}).
\par {\bf 2.2. Lemma.} {\it (1). If $G\in C^0(B_r\times X^{\otimes l},
L_{k-l}(X^{\otimes (k-l)};X))$, $\xi _i\in C^0(B_r,X)$
for each $i=1,...,k$ and $A_{l+i}\in C^0(B_r,
L(X))$ for each $i=1,...,k-l$, then
${\hat P}_{(\xi _{l+1},..,\xi _k)}[G(s;\xi _1,...,\xi _l)\circ
(A_{l+1}\otimes ...\otimes A_k)](v)\in C^0(B_r\times C^0(B_r,X)^{\otimes l},
C^0(B_r,X))$ 
as the function by $v, \xi _1,...,\xi _l$ for each fixed
$\xi _{l+1},..,\xi _k$ and $\hat P$ is of class $C^{\infty }$ by 
$\xi _{l+1},...,\xi _k$. 
\par (2). Moreover, if $G$ is of class of smoothness $C^m$ 
by arguments $\xi _1,...,\xi _l$, then ${\hat P}_{
(\xi _{l+1},....,\xi _k)}G$ is also in class of smoothness $C^m$
by $\xi _1,...,\xi _l$.}
\par {\bf Proof.} Since $B_r$ is compact, then $\xi _i$ are uniformly 
continuous together with $A_{l+i}(v)[\xi _{l+i}(v)].$
There is the following inequality
$$|\hat P_{(\eta _{l+1},...,\eta _k)}[G(v;\eta _1,...,\eta _l)\circ (
A_{l+1}\otimes ...\otimes A_k)](x)-$$
$$\hat P_{(\xi _{l+1},...,\xi _k)}[G(v;\xi _1,...,\xi _l)\circ (
A_{l+1}\otimes ...\otimes A_k)](y)|$$
$$\le \max (
|\hat P_{(\eta _{l+1},...,\eta _k)}[G(v;\eta _1,...,\eta _l)\circ (
A_{l+1}\otimes ...\otimes A_k)](x)-$$
$$\hat P_{(\xi _{l+1},...,\xi _k)}[G(v;\eta _1,...,\eta _l)\circ (
A_{l+1}\otimes ...\otimes A_k)](x)|,$$
$$|\hat P_{(\xi _{l+1},...,\xi _k)}[G(v;\eta _1,...,\eta _l)\circ (
A_{l+1}\otimes ...\otimes A_k)](x)-$$
$$\hat P_{(\xi _{l+1},...,\xi _k)}[G(v;\xi _1,...,\xi _l)\circ (
A_{l+1}\otimes ...\otimes A_k)](x)|,$$
$$|\hat P_{(\xi _{l+1},...,\xi _k)}[G(v;\xi _1,...,\xi _l)]\circ (
A_{l+1}\otimes ...\otimes A_k)](x)-$$
$$\hat P_{(\xi _{l+1},...,\xi _k)}[G(v;\xi _1,...,\xi _l)]\circ (
A_{l+1}\otimes ...\otimes A_k)](y)|.$$
In addition $\hat P$ is the linear operator by $\xi _{l+1},..., \xi _k$.
From this and Conditions $2.1.(i-iv)$ the first statement follows. 
The last statement follows from 
the linearity of $\hat P$ by $G$ and applying the operator of 
difference quotients $\bar \Phi ^m$ by $\xi _1,...,\xi _l$
(see \cite{luseamb,lubp2}).
\par {\bf 2.3. Lemma.} {\it If $\xi _i\in C^1(B_r,X)$ 
for each $i=1,...,k$
and Conditions (1) of Lemma 2.2 are satisfied, then
$${\hat P}_{(\xi _{l+1},...,\xi _k)}[G(s;\xi _1,...,\xi _l)\circ
(A_{l+1}\otimes ...\otimes A_k)](x) \in C^1(B_r,X)$$
as a function by the argument $x\in B_r$ and 
$$\partial /\partial x({\hat P}_{(\xi _{l+1},...,\xi _k)}
[G(s;\xi _1,...,\xi _l)\circ (A_{l+1}\otimes ...\otimes A_k)](x)=$$
$$\sum_{q=l+1}^k{\hat P}_{(\xi _{l+1},...,\xi _{q-1},\xi _{q+1},
...,\xi _k)}G(x;\xi _1,...,\xi _l).
(A_{l+1}(x){\xi }_{l+1}(x),...,A_{q-1}(x){\xi }_{q-1}(x),$$
$$A_q(x){\xi '}_q(x),A_{q+1}(x){\xi }_{q+1}(x),
...,A_k(x){\xi }_k(x))$$ such that
$$\| {\hat P}_{(\xi _{l+1},...,\xi _k)}[G(s;\xi _1,...,\xi _l)\circ
(A_{l+1}\otimes ...\otimes A_k)](x) \| _{C^1(B_r,X)}\le $$
$$\| G\| _{C^0(B_r\times X^{\otimes l},L_{k-l}(X^{\otimes (k-l)};X))}
\prod_{i=l+1}^k [\| A_i \| _{C^0(B_r,L(X))}
\| \xi _i\|_{C^1(B_r,X)}].$$}
\par {\bf Proof.} Let $\gamma :=
{\hat P}_{(\xi _{l+1},...,\xi _k)}[G(z;\xi _1,...,\xi _l)\circ
(A_{l+1}\otimes ...\otimes A_k)](x) - \mbox{ }
{\hat P}_{(\xi _{l+1},...,\xi _k)}[G(z;$ 
$\xi _1,...,\xi _l)\circ
(A_{l+1}\otimes ...\otimes A_k)](y) - \mbox{ }
(x-y)\sum_{q=l+1}^k 
{\hat P}_{(\xi _{l+1},...,\xi _{q-1},\xi _{q+1},
...,\xi _k)}[G(y;\xi _1,...,\xi _l).($ \\ $A_{l+1}(y)
\xi _{l+1}(y),...,A_{q-1}(y){\xi }_{q-1}(y),
A_q(y){\xi '}_q(y),A_{q+1}(y){\xi }_{q+1}(y),...,
A_k(y)\xi _k(y))]$ and \\
$\rho ^{s+1}\le |x-y| <\rho ^s$, where $s\in \bf N$.
Therefore, $x_0=y_0$,...,$x_s=y_s$, $x_{s+1}\ne y_{s+1}$ and 
$$\gamma =[\sum_{q=l+1}^kE(x_s)(v_{l+1},...,v_{q-1},h_q,z_{q+1},...,z_k)]
+E(x_s)(h_{l+1},h_{l+2},z_{l+3},...,z_k)$$
$$+E(x_s)(h_{l+1},v_{l+2},h_{l+3},z_{l+4},...,z_k)+...
+E(x_s)(v_{l+1},...,v_{k-2},h_{k-1},h_k)+...$$
$$+E(x_s)(h_{l+1},...,h_k) + \sum_{j=s+1}^{\infty }
\{ E(x_j)(\xi _{l+1}(x_{j+1})-\xi _{l+1}(x_j)),...,
(\xi _k(x_{j+1})-\xi _k(x_j)))$$
$$- E(y_j)(\xi _{l+1}(y_{j+1})-\xi _{l+1}(y_j)),...,
(\xi _k(y_{j+1})-\xi _k(y_j)))$$
$$-(x-y)\sum_{q=l+1}^k {\hat P}_{(\xi _{l+1},...,\xi _{q-1},\xi _{q+1},
...,\xi _k)}E(y)(\xi _{l+1}(y), ...,{\xi }_{q-1}(y),
{\xi '}_q(y),{\xi }_{q+1}(y),...,\xi _k(y)) ,$$
where
$v_j=\xi _j(x_{s+1})-\xi _j(x_s)$, $h_j=\xi _j(x_{s+1})-\xi _j(y_{s+1}),$
$z_j=\xi _j(y_{s+1})-\xi _j(y_s)$ for $j=l+1,...,k$ and
$$(i)\mbox{ }E:=E(x):=G(x;\xi _1,...,\xi _l).
(A_{l+1}(x)\otimes ...\otimes A_k(x))\mbox{ and}$$
$$(ii)\mbox{ }E(x)(\xi _{l+1},...,\xi _k):=G(x;\xi _1,...,\xi _l).
(A_{l+1}(x)\xi _{l+1},...,A_k(x)\xi _k)$$
in accordance with Formula 2.1.(3).
On the other hand, $\| \xi _i(y_{j+1})-\xi _i(y_j)
-(y_{j+1}-y_j)\xi _i(y)\| =\| (y_{j+1}-y_j)[(
{\bar \Phi }^1\xi _i)(y_j;1;y_{j+1}-y_j)- \xi _i(y)]\|
\le |y_{j+1}-y_j| \| \xi _i\| _{C^1(B_r,X)}$
and $E(x).(a_{l+1}+b_{l+1},...,a_k+b_k)-E(y).(a_{l+1},...,a_k)=$ \\
$E(x).(a_{l+1}+b_{l+1},....,a_k+b_k)-E(x)(a_{l+1},...,a_k)+
[E(x)-E(y)].(a_{l+1},....,a_k)=$ \\
$E(x).(b_{l+1},a_{l+2},...,a_k)+....
+E(x).(a_{l+1},....,a_{k-1},b_k)+E(x).(b_{l+1},b_{l+2},a_{l+3},...,a_k)+...$ \\
$+E(x).(a_{l+1},...,a_{k-2},b_{k-1},b_k)+...+E(x).(b_{l+1},...,b_k)+
[E(x)-E(y)].(a_{l+1},....,a_k)$ \\
for each $a_{l+1},...,a_k, b_{l+1},...,b_k \in C^0(B_r,X)$, hence \\
$\| [\sum_{q=l+1}^k
E(x_s)(v_{l+1},...,v_{q-1},h_q,z_{q+1},...,z_k)]-$ \\
$(x-y)\sum_{q=l+1}^k {\hat P}_{(\xi _{l+1},...,\xi _{q-1},\xi _{q+1},
...,\xi _k)}E(y)(\xi _{l+1}(y), ...,{\xi }_{q-1}(y),
{\xi '}_q(y),{\xi }_{q+1}(y),...,\xi _k(y)) \| $ \\   $\le  
\| E \|_{C^0}\rho ^s \prod_{q=l+1}^k \| \xi _q \| _{C^1} \alpha (s)$ and \\
$\| E(x_j)(\xi _{l+1}(x_{j+1})-\xi _{l+1}(x_j)),...,
(\xi _k(x_{j+1})-\xi _k(x_j)))$ \\
$- E(y_j)(\xi _{l+1}(y_{j+1})-\xi _{l+1}(y_j)),...,$ 
$(\xi _k(y_{j+1})-\xi _k(y_j))) \| $ \\
$\le \| E \|_{C^0}\rho ^s \prod_{q=l+1}^k \| \xi _q \| _{C^1} \alpha (s)$
for each $j\ge s+1$, where $\lim_{s\to \infty }\alpha (s)=0,$
consequently, $\lim_{|x-y|\to 0}\gamma =0$
and ${\bar \Phi }^1({\hat P}_{(\xi _{l+1},...,\xi _k)}E)(x)\in C^0(B_r,X)$, 
where $\Phi ^1\eta (x;h;\zeta )= \{ \eta (x+\zeta h)-
\eta (x) \} /\zeta $ for $0\ne \zeta \in \bf K$, $h\in H$, $\eta \in C^1
(U,Y)$, $U$ is open in $X$, $X$ and $Y$ are Banach spaces
over $\bf K$, ${\bar \Phi }^1\eta $ 
is a continuous extension of $\Phi ^1\eta $ on $U\times V\times 
B({\bf K},0,1)$ for a neighbourhood $V$ of $0$ in $X$ 
(see \S 2.3 \cite{luseamb} or \S I.2.3 \cite{lubp2}). Then 
$$(iii)\mbox{ }({\hat P}_{(\xi _{l+1},...,\xi _k)}E)(x)=
\sum_{n=0}^{\infty }(x_{n+1}-x_n)^{k-l}G(x_n;\xi _1,...,\xi _l).
(A_{l+1}(x_n)({\bar \Phi }^1\xi _{l+1})(x_n;$$
$$1;x_{n+1}-x_n),...,
(A_k(x_n)({\bar \Phi }^1\xi _k)(x_n;1;x_{n+1}-x_n)).$$
Let $\eta :=({\hat P}_wE)(x)-
({\hat P}_wE)(y)$, then
$\eta =E(x_s)(w(x_{s+1})-w(y_{s+1}))+\sum_{n=s+1}^{\infty }
\{ E(x_n)(w(x_{n+1})-w(x_n))-E(y_n)(w(y_{n+1})-w(y_n)) \} ,$
consequently, $\| \eta \| \le \| E\| _{C^0(B_r\times X^{\otimes l},L_{k-l}
(X^{\otimes (k-l)};X))}
(\prod_{i=l+1}^k \| \xi _i\| _{C^1}|x-y|)$, since
$E$ are polylinear mappings by $\xi _{l+1}(z),...,\xi _k(z)\in X$,
$|x_{s+1}-y_{s+1}|\le |x-y|$ and
$|x_{n+1}-x_n|\le |x-y|$ and
$|y_{n+1}-y_n |\le |x-y|$ for each $n>s$, where
$\rho ^{s+1}\le |x-y|<\rho ^s$, $w=(\xi _{l+1},...,\xi _k)$.
\par {\bf Note.} In particular, when $X=\bf K$, $l=0$, $k=1$, $A_1=1$ and 
$\xi (x)=x$ this gives the usual formula $d[\hat P_sG(s)](x)/dx=G(x).$
\par {\bf 2.4.} Suppose that $X$ and $Y$ are Banach spaces 
over a (complete relative to its uniformity) 
local field $\bf K$.
Let $X$ and $Y$ be isomorphic with the Banach spaces
$c_0(\alpha ,{\bf K})$ and $c_0(\beta ,{\bf K})$ 
and there are given the standard orthonormal bases $\{ e_j: j\in \alpha \} $
in $X$ and $\{ q_j: j\in \beta \} $ in $Y$ respectively, 
then each $E\in L(X,Y)$
has its matrix realisation $E_{j,k}:=q_k^*Ee_j$, where
$\alpha $ and $\beta $ are ordinals, $q_k^*\in Y^*$ is a continuous
$\bf K$-linear functional $q_k^*: Y\to \bf K$ corresponding to $q_k$
under the natural embedding $Y\hookrightarrow Y^*$
associated with the chosen basis, $Y^*$ is a topologically
conjugated or dual space of $\bf K$-linear functionals on $Y$.
\par {\bf 2.5.} Let $A$ be a commutative Banach algebra and $A^+$ denotes
the Gelfand space of $A$, that is, $A^+=Sp(A)$, where $Sp(A)$ 
in another words spectrum of $A$ was defined in
Chapter 6 \cite{roo}. Let $C_{\infty }(A^+,{\bf K})$ 
be the same space as in \cite{roo}. This means the following.
For a locally compact Hausdorff totally disconnected topological space
$E$ the vector space $C_{\infty }(E,{\bf K})$
is a subspace of a space $C(E,{\bf K})$
of bounded continuous functions $f: E\to \bf K$
such that for each $\epsilon >0$ there exists a compact subset
$V\subset  E$ for which $|f(x)|<\epsilon $ for each $x\in E\setminus V$.
When $E$ is not locally compact and have an embedding
into $B({\bf K},0,1)^{\gamma }$ (for example, when $\bf K$ is not 
locally compact) such that $E \cup \{ x_0 \}=cl(E)$
we put $C_{\infty }(E,{\bf K}):=\{ f\in C(E,{\bf K}): 
\lim_{x\to x_0}f(x)=0 \}$,
where $B(X,x,r):=\{ y\in X: d(x,y)\le r \} $ is a ball in the metric space
$(X,d)$, the closure $cl(E)$ is taken in $B({\bf K},0,1)^{\gamma }$, 
$\gamma $ is an ordinal, $x_0\in B({\bf K},0,1)^{\gamma }$.
\par {\bf Definition} (see also Ch. 6 in \cite{roo}). A commutative 
Banach algebra $A$ is called
a $C$-algebra if it is isomorphic with $C_{\infty }(X,{\bf K})$
for a locally compact Hausdorff totally disconnected topological
space $X$, where $f+g$ and $fg$
are defined pointwise for each $f, g\in C_{\infty }(X,{\bf K})$.
\par {\bf 2.6.} Let $H=c_0(\alpha ,{\bf K})$ and $X$ be a topological
space with the small inductive dimension $ind (X)=0$, where $\bf K$
is a complete field as the uniform space. 
A strong operator topology in $L(H,Y)$ (see \S 2.1)
is given by a base $V_{\epsilon ;E;x_1,...,x_n}:=\{ Z\in L(H,Y):
\sup_{1\le j\le n} \| (E-Z)x_j \|_Y <\epsilon \}$, where 
$0<\epsilon $, $E\in L(H,Y)$, $x_j\in H$; $j=1,...,n;$ $n\in \bf N$.
An $H$-projection-valued measure on an algebra
$\sf L$ of subsets of $X$ 
is a function $P$ on $\sf L$ assigning to each $A\in \sf L$
a projection $P(A)$ on $H$ and satisfying the following conditions:
\par $(i)$ $P(X)={\bf 1}_H$,
\par $(ii)$ for each sequence $\{ A_n : n\in {\bf N} \}$
of pairwise disjoint sets in $\sf L$ there 
are pairwise orthogonal projections $P(A_n)$
and $P(\bigcup_{n=1}^{\infty }A_n)=\sum_{n=1}^{\infty }P(A_n)$,
where ${\sf L}\supset Bco(X)$, $Bco(X)$ is an albegra of 
clopen (closed and open at the same time) subsets of $X$,
the convergence on the right hand side is unconditional in the
strong operator topology and the sum is equal to the projection
onto the closed linear span of $\bigcup_n\{ range(P(A_n)): n\in {\bf N} \} $
such that $P(\emptyset )=0$.
If $\eta \in H^*$ and $\xi \in H$, then $A\mapsto \eta (P(A)\xi )$
is a $\bf K$-valued measure on $\sf L$.
 Then by the definition $P(A)\le P(B)$ if and only if $range (P(A))
\subset range P(B)$. There are many projection operators on $H$,
but for $P$ there is chosen some such fixed system.
\par A subset $A\subset X$ is called $P$-null 
if there exists $B\in \sf L$ such that $A\subset B$ and $P(B)=0$,
$A$ is called $P$-measurable if $A\bigtriangleup B$ is $P$-null, 
where $A\bigtriangleup B:=(A\setminus B )\cup (B\setminus A)$.
A function $f: X\to \bf K$ is called $P$-measurable, if
$f^{-1}(D)$ is $P$-measurable for each $D$ in the algebra 
$Bco({\bf K})$ of clopen subsets of $\bf K$.
It is essentially bounded, if there exists $k>0$ such that
$\{ x: |f(x)|>k \} $ is $P$-null, $\| f\| _{\infty }$
is by the definition the infimum of such $k$. Then ${\sf F}:=
span_{\bf K} \{ Ch_B: B\in {\sf L} \} $ is called the 
space of simple functions,
where $Ch_B$ denotes the characteristic function of $B$.  The 
completion of $\sf F$ relative to $\| *\| _{\infty }$ is the Banach algebra
$L_{\infty }(P)$ under the pointwise multiplication.
\par For each
$f\in L_{\infty }(P)$ there exists the unique linear mapping
${\sf I}: {\sf F}\to L(H)$ by the following formula:
\par $(iii)$ ${\sf I}(\sum_{i=1}^n\lambda _i Ch_{B_i})=\sum_{i=1}^n
\lambda _iP(B_i)$, where $n\in \bf N$, 
$B_i\in \sf L$, $\lambda _i\in \bf K$.
Since
\par $(iv)$ $\| {\sf I}(f) \| =\| f\| _{\infty }$, then
$\sf I$ extends to a linear isometry (also called $\sf I$) of 
$L_{\infty }(P)$ onto $L(H)$.
\par If $f\in L_{\infty }(P)$, then the operator ${\sf I}(f)$
in $L(H)$ is called the spectral integral of $f$ with respect
to $P$ and it is denoted by
\par $(v)$ $\int_X f(x)P(dx):={\sf I}(f)$. 
\par From this definition using Chapter 7 \cite{roo} we get 
the following statement (compare with the classical case 
\S II.11.8 \cite{fell}).
\par {\bf Proposition.} {\it $(I)$. $\int_Xf(x)P(dx)=\int_Xg(x)P(dx)$
if and only if $f$ and $g$ differ only on a $P$-null set.
\par $(II)$. $\int_Xf(x)P(dx)$ is linear in $f$.
\par $(III)$. $\int_Xf(x)g(x)P(dx)=(\int_Xf(x)P(dx))(\int_Xg(x)P(dx))$
for each $f$ and $g\in L_{\infty }(P)$.
\par $(V)$. $\| \int_Xf(x)P(dx)\| =\| f\|_{\infty }$.
\par $(VI)$. If $A\in \sf L$, then $\int_XCh_A(x)P(dx)=P(A)$,
\par in particular $\int_XP(dx)=P(X)={\bf 1}_H$.
\par $(VII).$ For each pair $\xi \in H$ and $\eta ^*\in H^*$, 
let $\mu _{\xi ,\eta }(A):=
\eta ^*(P(A)\xi )$ for each $A\in \sf L$. If $E=\int_Xf(x)P(dx)$
then $\eta ^*(E\xi )=\int_Xf(x)\mu_{\xi ,\eta }(dx)$.
\par $(VIII).$ If $A\in \sf L$, then $P(A)$ commutes with
$\int_Xf(x)P(dx)$.}
\par An $H$-projection-valued measure $P$ on $Bco(X)$ is called
an $H$-projection-valued measure on $X$. We call
$P$ regular if 
\par $(v)$ $P(A)=\sup \{ P(C): C\subset A\mbox{ and }C \mbox{ is
compact } \} $ for each $A\in Bco(X)$, where $\sup $ is the least
closed subspace of $H$ containing $range\mbox{ }P(C)$ and to it
corresponds projector on this subspace. Indeed, $P(A)H$ is closed in $H$,
since $P^2(A)=P(A)$. Therefore, 
\par $(vi)$ $P(A)=\inf \{ P(U): U\mbox{ is open and }
U\supset A \} =I- \sup \{ P(C): C\subset X\setminus A\mbox{ and }
C\mbox{ is compact } \} $, hence 
\par $(vii)$ the infimum corresponds to the
projection on $\bigcap_{U\supset A, U\mbox{ is open}}P(U)H$.
\par A measure $\mu : Bco(X)\to \bf K$ is called regular, if for each 
$\epsilon >0$ and each $A\in Bco(X)$ with $\| A\|_{\mu }<\infty $ there
exists a compact subset $C\subset A$ such that 
$\| A\setminus C\|_{\mu }<\epsilon $. Since $\| P(X)\| =1$, 
then $\| \mu _{\xi ,\eta }\| \le \| \xi \|_H \| \eta \|_{H^*}$.
For the space $H$ over $\bf K$
measures $\mu _{\xi ,\eta }$ on the locally compact Hausdorff
totally disconnected topological space $X$ are tight 
for each $\xi , \eta $ in a subset $J\subset H\hookrightarrow H^*$ 
separating points of $H$
if and only if $P$ is defined on $Bco(X)$; $P$ is regular if and only if
$\mu _{\xi ,\eta }$ are regular for each $\xi , \eta \in J$
due to Conditions $(vi)$ and $(vii)$.
We can restrict our consideration by $\mu _{\xi ,\xi }$ 
instead of $\mu _{\xi ,\eta }$ with $\xi , \eta \in
span_{\bf K}J$, since ${+\choose -}2\mu _{\xi ,\eta }
=\mu _{\xi {+\choose -}\eta ,\xi {+\choose -}\eta }-\mu _{\xi ,\xi }
-\mu _{\eta ,\eta }$.
\par By the closed support of an $H$-projection-valued measure
$P$ on $X$ we mean the closed set $D$ of all those $x\in X$ such that
$P(U)\ne 0$ for each open neighbourhood $x\in U$, $supp\mbox{ }(P):=D$.
\par {\bf 2.7. Remark.} We fix a locally compact totally 
disconnected Hausdorff space $X$
and a Banach space $H$ over $\bf K$ and let $T: C_{\infty }(X,{\bf K})
\to L(H)$ be a linear continuous map from the $C$-algebra 
$C_{\infty }(X,{\bf K})$ of functions $f: X\to \bf K$ such that:
\par $(i)$  $T_{fg}=T_fT_g$ for each $f$ and 
$g\in C_{\infty }(X,{\bf K})$,
\par $(ii)$ $T_{\bf 1}=I$.
\par From this definition it follows, that $\| T\| \le 1$, since
$T_{f^n}=T_f^n$  for each $n\in \bf Z$ and $f\in C_{\infty }(X,{\bf K})$.
If $X$ is not compact and it is locally compact, 
then $X_{\infty }:=X\cup \{ x_{\infty } \}$
be its one-point Alexandroff compactification. 
Each $f\in C(X_{\infty },{\bf K})$
can be written just in one way in the form $f=\lambda {\bf 1}+g$,
where $g\in C_{\infty }(X,{\bf K})$ and $\bf 1$ is the unit function
on $X_{\infty }$. Therefore, we can extend $T: C_{\infty }(X,{\bf K})
\to L(H)$
to a linear map $T': C(X_{\infty },{\bf K})\to L(H)$ by setting
${T'}_{\lambda {\bf 1}+g}=\lambda {\bf 1}_H+T_g$ such that 
${T'}_{\bf 1}={\bf 1}_H$.
\par Therefore, $f\mapsto \eta ^*(T_f\xi )=:\tilde \mu _{\xi ,\eta }(f)$ 
is a continuous $\bf K$-linear
functional on $C_{\infty }(X,{\bf K})$, where $\xi \in H$ and $\eta ^*
\in H^*$.
In view of the Theorems 7.18 and 7.22 \cite{roo} about correspondence
between measures and continuous linear functionals 
(the non-Archimedean analog of the F. Riesz 
representation theorem) there exists the unique measure
$\mu _{\xi ,\eta }\in M(X)$ such that 
\par $(I)$ $\eta ^*(T_f\xi )=\int_X f(x)\mu _{\xi ,\eta }(dx)$ for each
$f\in C_{\infty }(X,{\bf K})$. Since
$T_{\bf 1}=I$, then $\mu _{\xi ,\eta }(X)=\eta ^*(\xi )=\xi ^*(\eta )$.
Then for each $A\in Bco(X)$ we have $\| A\|_{\mu _{\xi ,\eta }}
\le \| \xi \| \mbox{ } \| \eta \| \mbox{ }\sup_{f\ne 0} \| T_f\|$
$\le \| \xi \| \mbox{ } \| \eta \| $. Since $H$ considered as a subspace
of $H^*$ separates points in 
$H$, then for each $A\in Bco(X)$ there exists the unique linear operator
$P(A)\in L(H)$ such that:
\par $(II)$  $\| P(A)\| \le 1$ and  $\eta ^*(P(A)\xi )=\mu _{\xi ,\eta }(A)$,
since $\mu _{\xi ,\eta }(A)$ is a continuous bilinear $\bf K$-valued
functional by $\xi $ and $\eta \in H$.
From the existence of a $H$-projection-valued 
measure in the case of compact $X$
we get a projection-valued measure $P'$ on $X_{\infty }$ such that
${T'}_f=\int_{X_{\infty }}f(x)P'(dx)$ for each 
$f\in C(X_{\infty },{\bf K})$.
Suppose further in the locally compact non-compact case of $X$, that
\par $(iv)$ $span_{\bf K} \{ T_h\xi : f\in C_{\infty }(X,{\bf K}), 
\xi \in H \} $
is dense in $H$. From this last condition it follows, that
\par $(III)$ $P={P'}|_{Bco(X)}$ (see also \cite{put,roo}).
\par {\bf 2.8. Note.} A particular case of $H=C_{\infty }(X,{\bf K})$
for locally compact totally disconnected Hausdorff space $X$
and $T_f=f$ for each $f\in C_{\infty }(X,{\bf K})$ 
can be considered independently. 
Each such $f$ is a limit of a certain sequence by $n\in \bf N$
of finite sums
$\sum_jf(x_{j,n})Ch_{V_{j,n}}(x)$, where $\{ V_{j,n}:
j\in \Lambda _n \} $ is a finite partition of $X$ into the disjoint union
of subsets $V_{j,n}$ clopen in $X$, $x_{j,n}\in V_{j,n}$, 
$\Lambda _n\subset {\bf N}$, since $Range\mbox{ }(f)$ is bounded.
If to take $P(V)=Ch_V$ for each $V\in Bf(X)$, then
$T_fg=\lim_{n\to \infty }\sum_jf(x_{j,n})Ch_{V_{j,n}}(x)g
=\int_Xf(x)P(dx)g$
for each $g\in H$, so there is the bijective correspondence between
elements $f\in \sf A$ of a $C$-algebra $\sf A$ realised as 
$C_{\infty }(X,{\bf K})$
with $X=Sp ({\sf A})$ and their spectral integral representations.
It can be lightly seen that $P(V_1\cap V_2)=Ch_{V_1\cap V_2}=Ch_{V_1}Ch_{V_2}$
$=P(V_1)P(V_2)=P(V_2)P(V_1)$ for each $V_j\in Bco(X)$. If $\{ V_j:
V_j\in Bco(X), j\in {\bf N} \} $ is a disjoint family, then $P(\bigcup_jV_j)g=Ch_{\bigcup_jV_j}g=
\sum_jCh_{V_j}g=$ $\sum_jP(V_j)g$ for each $g\in H$. Also $P(\emptyset )H=
Ch_{\emptyset }H=\{ 0 \} $ and $P(X)g=Ch_Xg=g$ for each $g\in H$.
Therefore, $P$ is indeed the $H$-projection-valued measure. 
\par Suppose now that $X$ is not locally compact, for example, $X=c_0(\omega _0,
{\bf S})$ with an infinite residue class field $k$ of 
a field $\bf S$. Then there are $f\in C_{\infty }(X,{\bf K})$ 
for which convergence of finite or even countable or of the cardinality
$card\mbox{ }(k)$ (which may be greater or equal to 
$card \mbox{ }({\bf R})$) sums $\sum_jf(x_{j,n})Ch_{V_{j,n}}$
becomes a problem for a disjoint family $\{ V_{j,n}: j \} $ of 
clopen in $X$ subsets, since $\| Ch_{V_{j,n}} \|_{C(X,{\bf K})}=1$
for each $j$ and $n$.
\par {\bf 2.9. Remark.} Fix a Banach space $H$ over a non-Archimedean complete 
field $\bf F$,
as above $L(H)$ denotes the Banach algebra of all bounded $\bf F$-linear
operators on $H$. If $b\in L(H)$ we write shortly $Sp(b)$ instead of
$Sp_{L(H)}(b):=cl(Sp(span_{\bf F}\{ b^n: n=1,2,3,... \} )) $ 
(see also \cite{roo}). 
\par It was proved in Theorem 2 \cite{put} in the case of $\bf F$ 
with the dicrete valuation group, that each continuous $\bf F$-linear
operator $A: E\to H$ with $\| A \| \le 1$ from one Banach space $E$
into another $H$ has the form
$$A=U\sum_{n=0}^{\infty }\pi ^nP_{n,A},$$
where $P_n:=P_{n,A}$, $\{ P_n: n\ge 0 \} $ is a family of projections
and $P_nP_m=0$ for each $n\ne m$, 
$\| P_n \| \le 1$ and $P_n^2=P_n$ for each $n$,
$U$ is a partially isometric operator, that is,
$U|_{cl(\sum_nP_n(E))}$ is isometric, $U|_{E\ominus
cl(\sum_nP_n(E))}=0$, $ker (U)\supset ker (A)$, $Im(U)=cl(Im(A))$,
$\pi \in \bf F$,
$|\pi |<1$ and $\pi $ is the generator of the valuation group of $\bf F$.
\par For $\bf F$ not necessarily with the discrete valuation group
and a completely continuous operator $A$
it was proved the Fredholm alternative for the operator $I+A$ \cite{grus}.
\par We restrict our attentation to the case of the local field $\bf F$,
consequently, $\bf F$ has the dicrete valuation group.
If $\| A \| >1$ we get 
$$(i)\quad A=\lambda _A U\sum_{n=0}^{\infty }\pi ^nP_{n,A},$$
where $\lambda _A\in \bf F$ and $|\lambda _A|=\| A \| $. 
In view of \S \S 2.6-2.8 this is the particular case
of the spectral integration on the disceret topological space $X$.
Evidently, for each $1\le r<\infty $ there exists $J\in L(H)$ for which
$$(ii)\quad 
\{ \sum_{n\ge 0}s_n^r dim_{\bf F}P_{n,J}(H)\} ^{1/r} <\infty $$
for $1\le r<\infty $, 
where $J$ has the spectral decomposition given by Formula $(i)$,
$s_n:=|\lambda _J| |\pi |^n \| P_n \| $.
Using this result it is possible to give the following definition.
\par {\bf 2.10.1.} {\bf Definition.} Let $E$ and $H$ be two normed $\bf
F$-linear spaces, where $\bf F$ is an infinite spherically complete field
with a nontrivial non-Archimedean valuation. The $\bf F$-linear operator
$A\in L(E,H)$ is called of class $L_q(E,H)$ if there exists $a_n\in E^*$ 
and $y_n\in H$ for each $n\in \bf N$ such that 
$$(i)\quad (\sum_{n=1}^{\infty }\| a_n \| ^q_{E^*}\| y_n \| _H^q)<\infty $$
and $A$ has the form
$$(ii)\quad Ax=\sum_{n=1}^{\infty }a_n(x)y_n$$ for each
$x\in E$, where $1\le q<\infty $. For each such $A$ we put
$$(iii)\quad \nu _q(A)=\inf \{ \sum_{n=1}^{\infty }\| a_n \|^q_{E^*}
\| y_n \| ^q_H \} ^{1/q},$$ where the infimum is taken by all such
representations $(ii)$ of $A$,
$$(iv)\quad \nu _{\infty }(A):=\| A \| $$ and
$L_{\infty }(E,H):=L(E,H)$.
\par {\bf 2.10.2. Proposition.} {\it $L_q(E,H)$ is the normed $\bf F$-linear
space with the norm $\nu _q$.}
\par {\bf Proof.} Let $A\in L_q(E,H)$ and $1\le q<\infty $,
since the case $q=\infty $ follows from its definition.
Then $A$ has the representation 2.10.1.(ii). Then due to the ultrametric
inequality 
$$\| A x \| _H\le \| x \| _E \sup_{n\in \bf N}
( \| a_n \| _{E^*}\| y_n \| _H)\le \| x \| _E (\sum_{n=1}^{\infty }\| a_n \|
^q_{E^*} \| y_n \| ^q_H)^{1/q},$$ 
hence $\sup_{x\ne 0}\| A x \| _H/ \| x \| _E=:\| A \| \le \nu _q(A)$.
\par Let now $A,S \in L_q(E,H)$, then there exists $0<\delta <\infty $ and two
representations $Ax=\sum_{n=1}^{\infty }a_n(x)y_n$ and  
$Sx=\sum_{m=1}^{\infty }b_m(x)z_m$ for which 
$$(\sum_{n=1}^{\infty }
\| a_n \| ^q_{E^*}\| y_n \| ^q_H)^{1/q}\le \nu _q(A)+\delta \mbox{ and} $$
$$(\sum_{n=1}^{\infty }
\| b_n \| ^q_{E^*}\| z_n \| ^q_H)^{1/q}\le \nu _q(S)+\delta ,\mbox{ hence} $$
$(A+S)x=\sum_{n=1}^{\infty }(a_n(x)y_n+b_n(x)z_n)$ and 
$$\nu _q(A+S)\le (\sum_{n=1}^{\infty }\| a_n \| ^q\| y_n \| ^q)^{1/q}
+(\sum_{n=1}^{\infty }\| b_n \| ^q\| z_n \| ^q)^{1/q}\le \nu _q(A)+\nu _q(S)+
2\delta $$
 due to the H\"older inequality.
\par {\bf 2.11.} {\bf Proposition.} {\it If $J\in L_q(H)$, $S\in L_r(H)$
are commuting operators, the field $\bf F$ is with the discrete
valuation group and $1/q+1/r=1/v$, then $JS\in L_v(H)$,
where $1\le q, r, v \le \infty $.}
\par {\bf Proof.} Since $\bf F$ is with the discerete valuation, then $J$ and
$S$ have the decompositions given by Formula 2.9.(i). Certainly each projector
$P_{n,J}$ and $P_{m,S}$ belongs to $L_1(H)$ and have the decomposition
given by Formula 2.10.1.(ii). The $\bf F$-linear span of $\bigcup_{n,m}
range (P_{n,J}P_{m,S})$ is dense in $H$. In particular, for
each $x\in range (P_{n,J}P_{m,S})$ we have 
$J^kS^lx=\lambda _J^k\lambda ^l_S\pi ^{nk+ml}P_{n,J}P_{m,S}x.$ 
Applying \S 2.9 to commuting operators
$J^k$ and $S^l$ for each $k,l \in \bf N$
and using the base of $H$ we get projectors $P_{n,J}$ and $P_{m,S}$ 
which commute for each $n$ and $m$, consequently, 
$JS=U_JU_S\lambda _J\lambda _S\sum_{n\ge 0, m\ge 0}
\pi ^{n+m}P_{n,J}P_{m,S}$, hence 
$U_{JS}=U_JU_S$, $\lambda _{JS}=\lambda _J\lambda _S$,
$P_{l,JS}=\sum_{n+m=l}P_{n,J}P_{m,S}$. In view of the 
H\"older inequality $\nu _v(JS)=\inf (\sum_{n=0}^{\infty }s_{n,JS}^v 
dim_{\bf F}P_{n,JS}(H))^{1/v}
\le \nu _q(J)\nu _r(S)$ (see \S IX.4 \cite{reed}).
\par {\bf 2.12.1.} {\bf Proposition.} {\it If $E$ is the normed space and
$H$ is the Banach space over the field $\bf F$ 
(complete relative to its uniformity),
then $L_r(E,H)$ is the Banach space
such that if $J, S\in L_r(E,H)$, 
then $$\| J+S\|_r \le \| J\|_r + \| S\|_r;$$
$$\| b J\| _r=|b|\mbox{ }\| J \|_r$$ for each $b\in \bf K;$
$\| J\|_r=0$ if and only if $J=0$,
where $1\le r\le \infty $, $\| * \| _q:=\nu _q(*)$.}
\par {\bf Proof.} In view of Proposition 2.10.2 it remains to prove
that $L_r(E,H)$ is complete, when $H$ is complete.
Let $\{ T_{\alpha } \} $ be a Cauchy net in $L_r(E,H)$, then
there exists $T\in L(E,H)$ such that $\lim_{\alpha }T_{\alpha }x=Tx$ for each
$x\in E$, since $L_r(E,H)\subset L(E,H)$ and $L(E,H)$ is complete.
We demonstrate that $T\in L_r(E,H)$ and $T_{\alpha }$ convereges to
$T$ relative to $\nu _r$ for $1\le r<\infty $. 
Let $\alpha _k$ be a monotone subsequence
in $ \{ \alpha \} $ such that $\nu _r^r(T_{\alpha }-T_{\beta })<2^{-k-2}$
for each $\alpha , \beta \ge \alpha _k$, where $k\in \bf N$.
Since $T_{\alpha _{k+1}}-T_{\alpha _k}\in L_r(E,H)$, then
$(T_{\alpha _{k+1}}-T_{\alpha _k})x=\sum_{n=1}^{\infty }a_{n,k}(x)y_{n,k}$ with
$\sum_{n=1}^{\infty }\| a_{n,k}\| ^r \| y_{n,k} \| ^r<2^{-k-2}$. Therefore,
$(T_{\alpha _{k+p}}-T_{\alpha _k})x=\sum_{h=k}^{k+p-1}
\sum_{n=1}^{\infty }a_{n,h}(x)y_{n,h}$ for each $p\in \bf N$, consequently,
using convergence while $p$ tends to $\infty $ 
we get $(T-T_{\alpha _k})x=\sum_{h=k}^{\infty }
\sum_{n=1}^{\infty }a_{n,h}(x)y_{n,h}.$ Then $\nu _r^r(T-T_{\alpha _k})\le
\sum_{h=k}^{\infty }\sum_{n=1}^{\infty }\| a_{n,h}\| ^r \|
y_{n,h}\| ^r\le 2^{-k-1},$ hence $T-T_{\alpha _k}\in
L_r(E,H)$ and inevitably $T\in L_r(E,H)$. Moreover,
$\nu _r(T-T_{\alpha })\le \nu _r(T-T_{\alpha _k})+
\nu _r(T_{\alpha _k}-T_{\alpha })\le 2^{-(k-1)/r}2$ 
for each $\alpha \ge \alpha _k$.
\par {\bf 2.12.2. Proposition.} {\it Let $E, H, G$ be normed spaces
over spherically complete $\bf F$. 
If $T\in L(E,H)$ and $S\in L_r(H,G)$, then $ST\in
L_r(E,G)$ and $\nu _r(ST)\le \nu _r(S)\| T \| $. If $T\in L_r(E,H)$ and $S\in
L(H,G)$, then $ST\in L_r(E,G)$ and $\nu _r(ST)\le \| S \| \nu _r(T).$}
\par {\bf Proof.} For each $\delta >0$ there are $b_n\in H^*$
and $z_n\in G$ such that $Sy=\sum_{n=1}^{\infty }b_n(y)z_n$ for each $y\in H$
and $\sum_{n=1}^{\infty }\| b_n \| ^r \| z_n \| ^r
\le \nu _r^r(S)+\delta .$ Therefore, $STx=\sum_{n=1}^{\infty }
T^*b_n(x)z_n$ for each $x\in E$, hence $\nu _r(ST)\le \sum_{n=1}^{\infty }
\| T^*b_n \| ^r \| z_n \| ^r\le \| T \| [\nu _r^r(S)+\delta ],$ since
$ \| T^*b_n(x)\| = |b_n(Tx)|\le \| b_n \|
\| Tx \| \le \| b_n \| \| T \| \| x \| $, where $T^*\in L(H^*,E^*)$ is the
adjoint operator such that $b(Tx)=:(T^*b)(x)$ for each $b\in H^*$ and $x\in E$.
The operator $T^*$ exists due to the Hahn-Banach theorem for normed spaces
over the spherically complete field $\bf F$ \cite{roo}. 
\par {\bf 2.12.3. Proposition.} {\it If $T\in L_r(E,H)$, then
$T^*\in L_r(H^*,E^*)$ and $\nu _r(T^*)\le \nu _r(T)$, where $E$ and $H$ are
over the spherically complete field $\bf F$.}
\par {\bf Proof.} For each $\delta >0$ there are $a_n\in E^*$ and
$y_n\in H$ such that $Tx=\sum_{n=1}^{\infty }a_n(x)y_n$ for each $x\in E$ and
$ \sum_{n=1}^{\infty } \| a_n \| ^r \| y_n \| ^r\le \nu _r^r(T)+\delta $.
Since $(T^*b)(x)=b(Tx)=\sum_{n=1}^{\infty }a_n(x)b(y_n)$
for each $b\in H^*$ and $x\in E$, then $T^*b=\sum_{n=1}^{\infty }
y_n^*(b)a_n,$ where $y_n^*(b):=b(y_n)$, that is correct due to
the Hahn-Banach theorem for $E$ and $H$ over the spherically complete field
$\bf F$ \cite{roo}. Therefore, $\nu ^r_r(T^*)\le \sum_{n=1}^{\infty }\| y_n \|
^r \| a_n \| ^r\le \nu _r^r(T)+\delta ,$ since $ \| y^* \| _{H^*}= \| y \| _H$ 
for each $y\in H$.
\par {\bf 2.13.} For a space $L_k(H_1,...,H_k;H)$ of $k$-linear
mappings of $H_1\otimes ... \otimes H_k$ into $H$ we have its embedding
into $L(E,H)$, where $E$ is a normed space $H_1\otimes ... \otimes H_k$
in its maximum norm topology for normed spaces $H_1,...,H_k, H$
over $\bf F$ (see \S \S 2.1, 2.10). Therefore,
we can define the following normed space \\
$L_{r,k}(H_1,...,H_k;H):=L_k(H_1,...,H_k;H)\cap L_r(E;H)$ in particular \\
$L_{r,k}(H^{\otimes k};H):=L_k(H^{\otimes k};H)\cap L_r(H^{\otimes k};H)$ and\\
$L_{\infty ,k}(H_1,...,H_k;H):=L_k(H_1,...,H_k;H)$
with the norm $\nu _r(J)=:\| J\|_r$, where $1\le r\le \infty .$
Certainly, $L_{r,k}\subset L_{q,k}$ for each $1\le r<q\le \infty $.
\par Suppose that $(\Omega ,{\sf B},\lambda )$ is a probability space
(with non-negative measure $\lambda $), where $\sf B$
is a $\sigma $-algebra of subsets of $\Omega $.
We define a $\bf K$-linear Banach space $L^q(\Omega ,{\sf B},
\lambda ;L_{r,k}(H_1,...,H_k;H))$
and $L^q(\Omega ,{\sf B},\lambda ;L_k(H_1,...,H_k;H))$
as a completion of a family of mappings $\sum_{j=1}^nA_jCh_{W_j}$
with $A_j\in L_{r,k}(H_1,...,H_k;H)$ 
or $A_j\in L_k(H_1,...,H_k;H)$ respectively
and $W_j\in \sf B$ and $n\in \bf N$.
That is, 
as consisting of those mappings $\Omega \ni \nu \mapsto
A(\nu )\in L_{r,k}(H_1,...,H_k;H)$ for which $\| A(\nu )\|_r$
is $\lambda $-measurable and
$$\| A\|_{L^q}:=\{ \int_{\Omega }\| A(\nu )\|_r^q\lambda (d\nu ) \} ^{1/q}
<\infty ,$$ 
where $1\le q<\infty $;
$$\| A\|_{L^{\infty }}:=ess-\sup_{\lambda }\| A(\nu )\|_r.$$
\par {\bf 2.14.} We consider a $C^{\infty }$-manifold $X$ with an atlas
$At(X)=\{ (U_j,\phi _j): j \in \Lambda _X \} ,$ where $\bigcup_jU_j=X$,
$\phi _j(U_j)$ are open in $c_0(\alpha ,{\bf K})$ and $U_j$ are open in $X$,
$\phi _j : U_j\to \phi _j(U_j)$ are homeomorphisms, $\phi _i\circ \phi _j^{-1}
\in C^{\infty }$ for each $U_i\cap U_j\ne \emptyset $ and $\| \phi _i
\circ \phi _j^{-1} \| _{C^m}<\infty $ for each $m\in \bf N$,
$\phi _j(U_j)$ are bounded in $c_0(\alpha ,{\bf K})$ for each $j\in \Lambda _X$,
$\Lambda _X$ is a set, $C^n_b(X,H)$ is a completion of a set of all functions
$f: X\to H$ such that $f\circ \phi _j^{-1}\in C^n(\phi _j(U_j),H)$ for each
$j\in \Lambda _X$ and $\sup_j \| f\circ \phi _j ^{-1} \| _{C^n}=:\| f \|
_{C^n(X,H)}<\infty ,$ where $H$ is a Banach space over $\bf K$. Then
$C^n(X,H)$ is the set of all functions $f: X\to H$ such that
for each $x\in X$ there exists a neighbourhood $x\in U\subset X$ for which
$f|_U\in C^n_b(U,H)$. 
\par By $L^s(\Omega ,{\sf B},\lambda ;C^n(X,H))$
we denote a completion of a space of simple functions
$\sum_{j=1}^n\xi _j(x)Ch_{W_j}(\nu )$ with $\xi _j(x)\in C^n(X,H)$, 
$W_j\in \sf B$ and $n\in \bf N$, relative to the following 
norm 
$$\| \xi \| _{L^s} := \{ \int_{\Omega } \| \xi (x,\nu ) \| ^s_{C^n(X,H)}
\lambda (d\nu ) \} ^{1/s} < \infty $$
for each $1\le s<\infty $ or
$$\| \xi \|_{L^{\infty }}:=ess-\sup_{\lambda }
\| \xi (x,\nu )\|_{C^n(X,H)}<\infty ,$$  
where $X$ is the $C^{\infty }$ Banach manifold on $c_0(\alpha ,{\bf K})$, 
$\| \xi (x,\nu )\|_{C^n(X,H)}$ 
is attached to $\xi $ as a function by $x\in X$ with parameter
$\nu \in \Omega $ such that
$\| \xi (x,\nu )\|_{C^n(X,H)}$ is a measurable function by $\nu $. 
\par {\bf Theorem.} {\it Let $G\in L^r(\Omega ,{\sf B},\lambda ;C^0(B_R
\times H^{\otimes l},   
L_{k-l}(H^{\otimes (k-l)};H))$, $\xi _1,...,\xi _k\in L^q(\Omega ,{\sf B},
\lambda ;C^0(B_R,H))$, $A_{l+i}\in C^0(B_R,L(H))$ for each $i=1,...,k-l$ (see
\S 2.1),
where $B_R=B({\bf K},0,R)$, $G=G(x; \xi _1,...,\xi _l;\nu )$, 
$\xi _i=\xi _i(x,\nu )$
with $x\in B_R$, $\nu \in \Omega $, $1/r+1/q=1/s$ with 
$1\le r, q, s\le \infty $. Then
$({\hat P}_{(\xi _{l+1},...,\xi _k)}G\circ (A_{l+1}\otimes ... \otimes
A_k)\in L^s(\Omega ,{\sf B},\lambda ; C^0(B_R,H))$.}
\par {\bf Proof.} In $L^q(\Omega ,{\sf F},\lambda ;C^0(B_R\times V,W))$ the
family of step functions $f(t,x,\omega )=\sum_{j=1}^n
Ch_{U_j}(\omega )f_j(t,x)$ is dense, where $f_j\in C^0(B_R\times V,W)$, $Ch_U$
is the characteristic function of $U\in \sf F$, $n\in \bf N$, $V$ and $W$ are
Banach spaces over $\bf K$, $t\in B_R,$ $x\in V$, $\omega \in \Omega $,
since $\lambda (\Omega )=1$ and $\lambda $ is nonnegative \cite{boui,feder}.
Each matrix element $F_{h,b}(x,\nu )$ is in $L^r(
\Omega ,{\sf B},\lambda ;C^0(B_R,{\bf K}))$ and $\xi _j\in
L^q(\Omega ,{\sf B},\lambda ;C^0(B_R,{\bf K}))$, where \\
$F(x,\nu ):=G(x;a_1,...,a_l;\nu ).(A_{l+1}a_{l+1}(x),...,A_ka_k(x)),$ \\
$h\in H^*$, $b\in H$, $F_{h,b}:=h(Fb),$ $a_i\in C^0(B_R,H)$ for each 
$i=1,...,k$.
Since $\| \xi _j(x,\nu )\|_{C^0(X,H)} \in L^q(\lambda )$,
$\| F_{a,b} (x,\nu )\|_{C^n(X,H)} \in L^r(\lambda )$, then
$F(x,\nu ).w(x,\nu )\in L^s(\Omega ,{\sf B},\lambda ;C^0(B_R,H))$, 
where $w=(\xi _1,...,\xi _k)$ (see \S IX.4 \cite{reed}).
The operator $\hat P_wF$ is linear by $w$ and $F$, hence it is defined
on simple functions. In view of Lemma 2.2
$$\| \hat P_wF(x, \nu )\|_H
\le \| F(x,\nu )\|_{C^0(B_R\times H^{\otimes l},L_{k-l}(H^{\otimes (k-l)};H))}$$
$$\prod_{i=l+1}^k[ \| A_i \| _{C^0(B_R,L(H))}\| \xi _i(x,\nu )\|_{C^0(B_R,H)}]$$ 
for $\lambda $-a.e.
$\nu \in \Omega $, hence $\| (\hat P_wF)(x,\nu )\|_{L^s}
\le \| G \|_{L^r}\prod_{i=l+1}^k[\| A_i \| _{C^0}\| \xi _i\|_{L^q}]$.
\par {\bf Corollary.} {\it If in suppositions of Theorem 2.14
$\xi _i\in L^q(\Omega ,{\sf B},\lambda ;C^1(B_R,H))$ for 
each $i=1,...,k$, then $(\hat P_wF)\in L^s
(\Omega ,{\sf B},\lambda ;C^1(B_R,H))$ and
$$(i)\mbox{ }\| (\hat P_wG.(A_{l+1}\otimes ... \otimes A_k)) 
\|_{L^s(\lambda ;C^1(B_R,H))}
\le \| G\|_{L^r(\lambda ;C^0(B_R\times H^{\otimes l},L_{k-l}(H^{\otimes
(k-l)};H)))}$$ 
$$\prod_{i=l+1}^k[\| A_i \| _{C^0(B_R,L(H))}
\|\xi _i\|_{L^q(\lambda ;C^1(B_R,H))}.$$}
\par {\bf Proof.} In view of Lemma 2.3 and Theorem 2.14
$$\| (\hat P_wF)(x,\nu )
\|_{C^1(B_R,H)}\le \| G(x;\xi _1,...,\xi _l;\nu )\|_{C^0(B_R
\times H^{\otimes l},L_{k-l}(H^{\otimes (k-l)},H))}$$ 
$$\prod_{i=l+1}^k[\| A_i \| _{C^0(B_R,L(H))}\| \xi _i(x,\nu )\|_{C^1(B_R,H)}]$$ 
for $\lambda $-almost each $\nu \in \Omega $.
From this Formula $(i)$ follows.
\section{Markov quasimeasures for
a non-Archimedean Banach space.}
\par {\bf 3.1. Remark.} Let $H=c_0(\alpha ,{\bf K})$ be a Banach space over
a local field $\bf K$ with an ordinal $\alpha $
and the standard orthonormal base $\{ e_j: j\in \alpha \} $,
$e_j=(0,...,0,1,0,...)$ with $1$ on the $j$-th place.
Let $\sf U^P$ be a cylindrical algebra generated by projections on
finite-dimensional over $\bf K$ subspaces $F$ in $H$ and Borel 
$\sigma $-algebras $Bf(F)$. Denote by $\sf U$ the minimal $\sigma $-algebra
$\sigma ({\sf U^P})$ generated by $\sf U^P$. When
$card(\alpha )\le \aleph _0$, then ${\sf U}=Bf(H)$, where $card (A) $
denotes the cardinality of a set $A$.
Each vector $x\in H$ is considered as continuous linear functional 
on $H$ by the formula $x(y)=\sum_jx^jy^j$ for each $y\in H$, 
so there is the natural embedding $H\hookrightarrow H^*=l^{\infty }
(\alpha ,{\bf K})$, where $x=\sum_jx^je_j$, $x^j\in \bf K$.
\par {\bf 3.2. Notes and definitions.} Let $T=B({\bf K},t_0,r)$ 
be a ball in the field $\bf K$ of radius $r>0$ and containing a point $t_0$
and $X_t=X$ be a locally $\bf K$-convex space
for each $t\in T$. Put $(\tilde X_T,
{\tilde {\sf U}}):=\prod_{t\in T}(X_t,{\sf U}_t)$ be a product of
measurable spaces, where ${\sf U}_t$ are $\sigma $-algebras 
of subsets of $X_t$, $\tilde {\sf U}$ is the $\sigma $-algebra of cylindrical
subsets of $\tilde X_T$ generated by projections $\tilde \pi _q: \tilde X_t\to
X^q$, $X^q:=\prod_{t\in q}X_t,$ $q\subset T$ is a finite subset of $T$
(see \S I.1.3 \cite{dalf}).
A function $P(t_1,x_1,t_2,A)$ with values in $\bf C$
for each $t_1\ne t_2\in T$, $x_1\in X_{t_1}$, $A\in {\sf U}_{t_2}$
is called a transition measure if it satisfies the following conditions:\\
$(i)\mbox{ the set function }\nu _{x_1,t_1,t_2,}(A)
:=P(t_1,x_1,t_2,A)\mbox{ is 
a }\sigma \mbox{-additive measure on }$ $(X_{t_2},{\sf U}_{t_2});$\\
$(ii)\mbox{ the function }\alpha _{t_1,t_2,A}(x_1)
:=P(t_1,x_1,t_2,A)\mbox{ of the variable } x_1$ $\mbox{ is  }
{\sf U}_{t_1}-\mbox{measurable};$ \\
$$(iii)\mbox{ }P(t_1,x_1,t_2,A)=\int_{X_s}
P(t_1,x_1,s,dy)P(s,y,t_2,A)\mbox{ for each }t_1\ne t_2\in T.$$
A transition measure  $P(t_1,x_1,t_2,A)$ is called normalised if
$$(iv)\mbox{ }P(t_1,x_1,t_2,X_{t_2})=1\mbox{ for each }t_1\ne t_2\in T.$$
For each set $q=(t_0,t_1,..,t_{n+1})$ of pairwise distinct points in $T$
there is defined a measure in $X^s:=\prod_{t\in s}X_t$ by the formula
$$(v)\mbox{ }\mu ^q_{x_0}(E)=\int_E\prod_{k=1}^{n+1}P(t_{k-1},x_{k-1},
t_k,dx_k), \mbox{ }E\in {\sf U}^s:=\prod_{t\in s}{\sf U}_t,$$
where $s=q\setminus \{ t_0 \}$, variables $x_1,...,x_{n+1}$ are such that
$(x_1,...,x_{n+1})\in E$, $x_0\in X_{t_0}$ is fixed.
\par Let $E=E_1\times X_{t_j}\times E_2$, where 
$E_1\in \prod_{i=1}^{j-1}{\sf U}_{t_i}$, $E_2\in \prod_{i=j+1}^{n+1}
{\sf U}_{t_i}$, then
\par $(vi)$ $\mu ^q_{x_0}(E)=\int_{E_1\times E_2}[\prod_{k=1}^{j-1}P(
t_{k-1},x_{k-1},t_k,dx_k)]\times $ 
$[\int_{X_{t_j}}P(t_{j-1},x_{j-1},t_j,dx_j)$
$\prod_{k=j+1}^{n+1}P(t_{k-1},x_{k-1},t_k,dx_k)]$ 
$=\mu ^r_{x_0}(E_1\times E_2),$ where $r=q\setminus \{ t_j \} .$
From Equation $(vi)$ it follows, that 
$$(vii)\mbox{ }[\mu ^q_{x_0}]^{\pi ^q_v}=\mu ^v_{x_0}$$ 
for each $v<q$ (that is, $v\subset q$), where $\pi ^q_v: X^s\to X^w$
is the natural projection, $s=q\setminus \{ t_0 \} ,$
$w=v\setminus \{ t_0 \} .$ If the transition measure 
$P(t,x_1,t_2,dx_2)$ is normalised and $F=E\times X_{t_{n+1}}$ with
$E\in {\sf U}^s$, then
$$(viii)\mbox{ }{\mu }^v_{x_0}(E):=\mu ^q_{x_0}(F)
=\int_E\prod_{j=1}^nP(t_{j-1},x_{j-1},t_j,dx_j),$$
where $q=(t_0,...,t_{n+1})$, $v=(t_0,...,t_n)$, points $t_0,...,t_n$ are
pairwise distinct in $T$.
If $\nu $ is the complex-valued measure on $(X,{\sf U})$, then
$\nu =\nu _1-\nu _2+i\nu _3-i\nu _4,$ where $\nu _j$ is a nonnegative measure
on $(X,{\sf U})$ for each $j=1,...,4$, $i:=(-1)^{1/2}\in \bf C$,
$\sf U$ is a $\sigma $-algebra of subsets of $X$.
By the definition $ \| \nu \| :=\sum_{j=1}^4\nu _j(X)$ and it is called 
the variation of the measure $\nu $ on $X$.
Therefore, due to Conditions $(iv,v,vii):$ $\{ \mu ^q_{x_0}; \pi ^q_v;
\Upsilon _T \} $ is the consistent family of measures, which induce
the quasimeasure $\tilde \mu _{x_0}$ on $(\tilde X_T, \tilde {\sf U})$
such that $\tilde \mu _{x_0}(\pi
_q^{-1}(E))=\mu ^q_{x_0}(E)$ for each $E\in {\sf U}^s$, where 
$\Upsilon _T$ is the family of all finite subsets $q$ in $T$
such that $t_0\in q\subset T$, $v\le q\in \Upsilon _T$, $\pi _q: \tilde U_T\to
X^s$ is the natural projection, $s=q\setminus \{ t_0 \} $.
\par The quasimeasures given by Equations $(i-v,vii)$
are called Markov quasimeasures.
\par {\bf 3.3.} {\bf Proposition. 1.} {\it If a normalized transition measure
$P$ satisfies the condition
$$(i)\mbox{ }C:=\sup_q[\sum_{k=1}^nln (\sup_x
\| \nu _{x,t_{k-1},t_k} \|)] < \infty ,$$ where $q=(t_0,t_1,...,t_n)$
with pairwise distinct points $t_0,..,t_n \in T$ and $n\in \bf N$,
then the Markov quasimeasure ${\tilde \mu }_{x_0}$ is bounded.}
\par {\bf 3.3.2. Proposition.} {\it If 
$$(ii)\mbox{ }C_x:=\sup_q[\sum_{k=1}^nln 
\| \nu _{x,t_{k-1},t_k} \|]=\infty $$ 
for each $x,$ where $q=(t_0,t_1,...,t_n)$
with pairwise distinct points $t_0,..,t_n \in T$ and $n\in \bf N$,
then the Markov quasimeasure ${\tilde \mu }_{x_0}$ has the unbounded variation
on each nonvoid set $E\in {\sf U}^s.$}
\par {\bf Proof.} (1). If $E\in \tilde {\sf U}$, 
then $E\in {\sf U}^s$ for some set
$q=(t_0,t_1,...,t_n)$ with pairwise distinct points
$t_0,...,t_n \in T$ and $n\in \bf N$ and $s=q\setminus \{ t_0 \} $,
consequently,
$|\mu ^q_{x_0}(E)|\le \prod_{k=1}^n
\sup_x\| \nu _{x,t_{k-1},t_k}\|$ 
$\le exp(C),$ since $t_k\in T$ for each $k=0,1,...,n$.
\par (2). For each $(t_1,t_2,x)$ there exists a compact set
$\delta (t_1,t_2,x)\in {\sf U}_{t_2}$ such that
$P(t_1,x_1,t_2,\delta (t_1,t_2,x))>1+\epsilon (t_1,t_2,x_1,x)$,
where $\epsilon (t_1,t_2,x)>0$. In view of Condition $(ii)$ 
for each $R>0$ and $x$ we choose
$q$ such that $\sum_{k=1}^n\epsilon (t_k,t_{k+1},x_1,x)>R$. 
For chosen $u\ne u_1\in T$ and $x\in X_u$ we represent the set $\delta (u,u_1,x)$
as a finite union of disjoint subsets $\gamma _{j_1}$ such that for each
$\gamma _{j_1}$ and $u_2\ne u_1$ there is a set $\delta _{j_1}$
satisfying $P(u_1,x_1,u_2,\delta _{j_1})\ge 1+\epsilon (u_1,u_2,x_1,x)$
for each $x\in \gamma _{j_1}.$
Then by induction $\delta _{j_1,...,j_n}=\bigcup_{j_{n+1}=1}^{m_{n+1}}
\gamma _{j_1,...,j_{n+1}}$ so that for $u_{n+2}\ne u_{n+1}\in T$
there is a set $\delta _{j_1,...,j_{n+1}}$ for which
$P(u_{n+1},x_{n+1},u_{n+2},\delta _{j_1,...,j_{n+1}})\ge 1+
\epsilon (u_{n+1},u_{n+2},x_{n+1},x)$ for each $x\in 
\gamma _{j_1,...,j_{n+1}}.$
Put $\Gamma ^{u,x_0}_{j_1,...,j_n}=$ $\{ x: x(u)=x_0, x(u_1)
\in \gamma _{j_1},...,x(u_n)\in \delta _{j_1,...,j_n},$
$x(u_{n+1})\in \gamma _{j_1,...,j_n} \} $
and $\Gamma ^{u,x_0}:=\bigcup_{j_1,...,j_n}
\Gamma ^{u,x_0}_{j_1,...,j_n}$. Then ${\tilde \mu }_{x_0}(\Gamma ^{u,x_0})
=\sum_{j_1,...,j_n}\int_{\delta _{j_1,...,j_n}}\int_{\gamma _{j_1,...,j_n}}
...\int_{\gamma _{j_1}}\prod_{k=1}^{n+1}P(u_{k-1},x_{k-1},u_k,dx_k)\ge  \\
\prod_{k=1}^n[1+\epsilon (u_{k-1},u_k,x_{k-1},x_k)]>R$.
\par {\bf 3.3.3} Evidently Condition $(i)$ of Proposition 3.3.1
is satisfied for the nonnegative normalized transition measure.
\par {\bf 3.4.} Let $X_t=X$ for each $t\in T$, ${\tilde X}_{t_0,x_0}:=
\{ x\in {\tilde X}_T: $ $x(t_0)=x_0 \} .$ We define a projection
operator $\bar \pi _q:$ $x\mapsto x_q$, where $x_q$ is defined on 
$q=(t_0,...,t_{n+1})$ such that $x_q(t)=x(t)$ for each $t\in q$,
that is, $x_q=x|_q$. For every $F: {\tilde X}_T\to \bf C$
there corresponds $(S_qF)(x):=F(x_q)=F_q(y_0,...,y_n),$ where $y_j=x(t_j)$.
$F_q: X^q\to \bf C$. We put ${\sf F}:=$ $\{ F| F: {\tilde X}_T\to {\bf C},$
$S_qF\mbox{ are }{\sf U}^q-\mbox{measurable} \} $.
If $F\in \sf F$, $\tau =t_0\in q$, then there exists an integral
$$(i)\mbox{ }J_q(F)=\int_{X^q}(S_qF)(x_0,...,x_n)\prod_{k=1}^{n+1}P(t_{k-1},
x_{k-1},t_k,dx_k).$$
\par {\bf Definition.} A function $F$ is called integrable with respect 
to the Markov quasimeasure $\mu _{x_0}$ if the limit
$$(ii)\mbox{ }\lim_qJ_q(F)=:J(F)$$ along the generalized net by 
finite subsets $q$ of $T$ exists. This limit is called a functional 
integral with respect to the Markov quasimeasure:
$$(iii)\mbox{ }J(F)=\int_{{\tilde X}_{t_0,x_0}}F(x)\mu _{x_0}(dx).$$
\par {\bf 3.5. Remark.} Consider a complex-valued measure 
$P(t,A)$ on $(X,{\sf U})$ for each $t\in T:=B({\bf K},0,R)$
such that $A-x\in \sf U$ for each $A\in \sf U$ and $x\in X$, 
where $A\in \sf U$, $X$ is a locally $\bf K$-convex space,
$\sf U$ is a $\sigma $-algebra of $X$.
Suppose $P$ be a spatially homogeneous transition measure (see also \S 3.2), 
that is,
$$(i)\mbox{ }P(t_1,x_1,t_2,A)=P(t_2-t_1,A-x_1)$$ 
for each $A\in \sf U$, $t_1\ne t_2 \in T$ and every $x_1\in X$, where $P(t,A)$ 
satisfies the following condition:
$$(ii)\mbox{ }P(t_1+t_2,A)=\int_XP(t_1,dy)P(t_2,A-y).$$
Such a transition measure $P(t_1,x_1,t_2,A)$ is 
called homogeneous. In particular for $T=\bf Z_p$
we have 
$$(iii)\mbox{ }P(t+1,A)=\int_XP(t,dy)P(1,A-y).$$ 
If $P(t,A)$ is a continuous function by $t\in T$ 
for each fixed $A\in \sf U$, then Equation $(iii)$ defines
$P(t,A)$ for each $t\in T$, when $P(1,A)$ is given, since
$\bf Z$ is dense in $\bf Z_p$.
\par {\bf 3.6. Notes and definition.}
 Let $X$ be a locally $\bf K$-convex space
and $P$ satisfies Conditions $3.2(i-iii)$. 
For $x$ and $z \in \bf Q_p^n$ we denote by
$(z,x)$ the following sum $\sum_{j=1}^n x_jz_j$, where $x=(x_j:$
$j=1,...,n)$, $x_j \in \bf Q_p$. Each number $y\in \bf Q_p$ has  a
decomposition $y=\sum_l
a_lp^l$, where $a_l \in (0,1,...,p-1)$,
$\min (l:$ $a_l\ne 0)=:ord_p(y)> - \infty $
($ord(0):=\infty $) \cite{nari,sch1}, 
we define a symbol $\{ y \}_p:=\sum_{l<0} a_lp^l$ for
$|y|_p>1$ and $\{ y\} _p=0$ for $|y|_p \le 1$. 
We consider a character 
of $X$, $\chi _{\gamma }: X\to {\bf C}$  given by the following formula:
$$(i)\mbox{ }\chi _{\gamma }(x)=\epsilon ^{z^{-1}\{ (e,\gamma (x)) \} _p}$$ 
for each $ \{ (e,\gamma (x)) \} _p\ne 0$,
$\chi _{\gamma }(x):=1$ for $ \{ (e,\gamma (x)) \} _p=0,$
where $\epsilon =1^z$ is a root of unity, $z=p^{ord(\{ (e,\gamma (x)) \} _p)},$ 
$\gamma \in X^*$, $X^*$ denotes the topologically 
conjugated space of continuous $\bf K$-linear functionals
on $X$, the field $\bf K$ as the $\bf Q_p$-linear space is
$n$-dimensional, that is,
$dim_{\bf Q_p}{\bf K}=n$, $\bf K$ as the Banach space over $\bf Q_p$
is isomorphic with $\bf Q_p^n$, $e=(1,...,1)\in \bf Q_p^n$
(see \cite{vla3} and \cite{lu6}). Then 
$$(ii)\mbox{ }\phi (t_1,x_1,t_2,y):=\int_X\chi _y(x)P(t_1,x_1,t_2,dx)$$
is the characteristic functional of the transition measure $P(t_1,x_1,t_2,dx)$
for each $t_1\ne t_2\in T=B({\bf K},t_0,R)$ and each $x_1\in X$.
In the particular case of $P$ satisfying Conditions $3.5.(i,ii)$ with $t_0=0$
its characteristic functional is such that 
$$(iii)\quad \phi (t_1,x_1,t_2,y)=\psi (t_2-t_1,y)\chi _y(x_1),\mbox{ where}$$
$$(iv)\quad \psi (t,y):=\int_X\chi _y(x)P(t,dx)\mbox{ and}$$
$$(v)\quad \psi (t_1+t_2,y)=\psi (t_1,y)\psi (t_2,y)$$ 
for each $t_1\ne t_2\in T$ and $y\in X^*$, $x_1\in X$.
\section{Non-Archimedean stochastic processes.}
\par {\bf 4.1. Remark and definition.} A measurable space $(\Omega ,{\sf F})$ 
with a normalised 
non-negative measure $\lambda $ on a $\sigma $-algebra $\sf F$ of a set
$\Omega $
is called a probability space and is denoted by
$(\Omega ,{\sf F},\lambda )$. Points $\omega \in \Omega $ are called 
elementary events and values $\lambda (S)$  
probabilities of events $S\in \sf F$. A measurable map
$\xi : (\Omega ,{\sf F})\to (X,{\sf B})$ is called a random variable
with values in $X$, where $\sf B$ is the $\sigma $-algebra
of a locally $\bf K$-convex space $X$. 
The random variable $\xi $ induces a normalized measure $\nu _{\xi }(A):=
\lambda (\xi ^{-1}(A))$ in $X$ and a new probability space 
$(X,{\sf B},\nu _{\xi }).$
We take $X=C^0(T,H)$ (see \S 2.1) and the $\sigma $-algebra
${\sf B}$ which is the subalgebra of the Borel $\sigma $-algebra $Bf(X)$ 
of $X$, where $H$ is a Banach space over $\bf K$, 
$T=B({\bf K},t_0,R)=:B_R,$ $0<R<\infty ,$ $\bf K$ is the local field.
A random variable
$\xi : \omega \mapsto \xi (t,\omega )$ with values in $(X,{\sf B})$
is called a (non-Archimedean) stochastic process on $T$ with values in $H$.
\par Events $S_1,...,S_n$ are called independent in total if
$P(\prod_{k=1}^nS_k)=\prod_{k=1}^nP(S_k)$. $\sigma $-Subalgebras
${\sf F}_k\subset {\sf F}$ are said to be independent if
all collections of events $S_k\in {\sf F}_k$ are independent in total, 
where $k=1,...,n$, $n\in \bf N$. To each collection of random variables
$\xi _{\gamma }$ on $(\Omega ,{\sf F})$ with $\gamma \in \Upsilon $
is related the minimal $\sigma $-algebra ${\sf F}_{\Upsilon }\subset \sf F$
with respect to which all $\xi _{\gamma }$ are measurable, where $\Upsilon $
is a set.
Collections $\{ \xi _{\gamma }: $ $\gamma \in \Upsilon _j \} $
are called independent if such are 
${\sf F}_{\Upsilon _j}$, where $\Upsilon _j\subset \Upsilon $ for each
$j=1,...,n,$ $n\in \bf N$.
\par {\bf 4.2. Defintion.} We define a (non-Archimedean)
stochastic process 
$w(t,\omega )$ with values in $H$ 
as a stochastic process such that:
\par $(i)$ the differences $w(t_4,\omega )-w(t_3,\omega )$ 
and $w(t_2,\omega )-w(t_1,\omega )$ are independent
for each chosen $\omega $, $(t_1,t_2)$ and $(t_3,t_4)$ with $t_1\ne t_2$,
$t_3\ne t_4$, either $t_1$ or $t_2$ is not in the two-element set
$ \{ t_3,t_4 \} ,$ where $\omega \in \Omega ;$
\par $(ii)$ the random variable $\omega (t,\omega )-\omega (u,\omega )$ has 
a distribution $\mu ^{F_{t,u}},$ where $\mu $ is a probability 
measure on $C^0(T,H)$, 
$\mu ^g(A):=\mu (g^{-1}(A))$ for $g\in C^0(T,H)^*$ and each $A\in
\sf B$, a continuous linear functional $F_{t,u}$ is given by
the formula $F_{t,u}(w):=w(t,\omega )-w(u,\omega )$ 
for each $w\in L^q(\Omega ,{\sf F},\lambda ;C^0_0(T,H)),$
where $1\le q\le \infty ,$ $C^0_0(T,H):=\{ f: f\in C^0(T,H), f(t_0)=0 \} $ 
is the closed subspace of $C^0_0(T,H)$.
\par $(iii)$ we also put $w(0,\omega )=0,$  
that is, we consider a Banach subspace
$L^q(\Omega ,{\sf F},\lambda ;C^0_0(T,H))$
of $L^q(\Omega ,{\sf F},\lambda ;C^0(T,H))$,
where $\Omega \ne \emptyset $.
\par This definiton is justified by the following Theorem.
\par {\bf 4.3.} {\bf Theorem.} {\it There exists a family 
of pairwise inequivalent (non-Archimedean) stochastic processes on $C^0_0(T,H)$
of the cardinality ${\sf c}$, where ${\sf c}:=card ({\bf R})$.}
\par {\bf Proof.} Since $H$ is over the local field, then $H$ has
a projection $\pi _0$ on its Banach subspace $H_0$ of separable type
over $\bf K$ (see its definition in \cite{roo}), that is,
$H_0$ is isomorphic with $c_0(\alpha ,{\bf K})$ with countable $\alpha $.
Therefore, a $\sigma $-additive measure $\mu _0$ on $(H_0,Bf(H_0))$
induces a $\sigma $-additive measure $\mu $ on $(H, \pi _0^{-1}[Bf(H_0)]),$
where $\pi _0^{-1}[Bf(H_0)]:=\{ \pi _0^{-1}(A): A\in Bf(H_0) \} .$
Therefore, it is sufficient to consider
the case of $H$ of separable type over $\bf K$.
\par If $w$ is the real-valued nonnegative Haar measure on $\bf K$
with $w(B({\bf K},0,1))=1$,
then it has not any atoms, since it is defined on
$Bf({\bf K})$, each singleton $\{ x\} $ is the Borel subset
and $w(y+A)=w(A)$ for each $A\in Bf({\bf K})$. Indeed,
if $w$ would have some atom $E$, then it would be a singleton,
since $\bf K$ is the complete separable metric space
and for each disjoint $w$-measurable subsets $A$ and $S$
in $E$ either $w(A)=w(E)>0$ with $w(S)=0$ or
$w(S)=w(E)>0$ with $w(A)=0$. But $\sum_{y\in {\bf K}}w(y+\{ x \} )=
\infty $, when $w(\{ x \} )>0$ for a singleton $\{ x \} $
(see Chapter VII in \cite{boui}).
Therefore, each measure 
$\mu _j(dx^j)=f_j(x^j)w(dx^j)$
on $\bf K$ has not any atom, since $w$ has not any atom, where
$f_j\in L^1({\bf K},Bf({\bf K}),w,{\bf R})$ (that is, $f_j$ is $w$-measurable
and $\| f _j \| _{L^1}:=\int_{\bf K}|f_j(x)|w(dx)<\infty $)
and $\mu _j({\bf K})=1$.
Hence each measure $\mu $ on $C^0_0(T,H)$ has not any atom, when
$\mu (dx)=\bigotimes_{j=1}^{\infty }\mu _j(dx^j),$ where 
$C^0_0(T,H)$ is isomorphic with
$c_0(\omega _0,{\bf K})$, $x\in C^0_0(T,H)$, $x=(x^j: j\in \omega _0),$ 
$x^j\in \bf K$, $x=\sum_jx^je_j$, $e_j$ is the standard othonormal base in 
$c_0(\omega _0,{\bf K})$, $\omega _0$ is the first countable ordinal,
since $\bf K$ is the local field (see \cite{roo} and
\cite{lu6}).
\par Let on the Banach space $c_0:=c_0(\omega _0,{\bf K})$ 
there is given an operator $J\in L_1(c_0)$ such that
$Je_i=v_ie_i$ with $v_i\ne 0$ for each $i$ and a measure $\nu (dx):=f(x)w(dx)$,
where $f: {\bf K}\to [0,1]$ is a function belonging to the space
$L^1({\bf K},w,{\bf R})$ such that $\lim_{|x|\to \infty }f(x)=0$
and $\nu ({\bf K})=1$, $\nu (S)>0$ for each open subset $S$ in $\bf K$,
for example, when $f(x)>0$ $w$-almost everywhere. 
In view of the Prohorov theorem 
there exists a $\sigma $-additive product
measure \\
$(i)\quad \mu (dx):=\prod_{i=1}^{\infty }\nu _i(dx^i)$ on the $\sigma $-algebra
of Borel subsets of $c_0$, since the Borel $\sigma $-algebras defined for the
weak topology of $c_0$ and for the norm topology of $c_0$ coincide,
where $\nu _i(dx^i):=f(x^i/v_i)\nu (dx^i/v_i)$ (see \cite{boui,lu6}).
\par Let $Z$ be a compact subset without
isolated points in a local field $\bf K$, for example,
$Z=B({\bf K},t_0,1)$. Then 
the Banach space $C^0(Z,{\bf K})$ has the Amice polynomial 
orthonormal base
$Q_m(x)$, where $x\in Z$, $m\in {\bf N_o}:=\{ 0,1,2,... \} $ \cite{ami}.
Each $f\in C^0$ has a decomposition $f(x)=\sum_ma_m(f)Q_m(x)$
such that $\lim_{m\to \infty }a_m=0$, where $a_m\in \bf K$.
These decompostions establish the isometric isomorphism
$\theta : C^0(T,{\bf K})\to c_0(\omega _0,{\bf K})$
such that $\| f\|_{C^0}=\max_m|a_m(f)|=\| \theta (f)\|_{c_0}$.
\par If $H=c_0(\omega _0,{\bf K}),$ then the Banach space
$C^0(T,H)$ is isomorphic with the tensor product 
$C^0(T,{\bf K})\otimes H$ (see \S 4.R \cite{roo}).
If $J_i\in L_1(Y_i)$ is nondegenerate for each $i=1,2$, that is, 
$ker (J_i)=\{ 0 \} $, then $J:=J_1\otimes J_2\in L_1(Y_1\otimes Y_2)$
is nondegenerate (see also Theorem 4.33 \cite{roo}).
If $u_i$ are roots of basic polynomils $Q_m$ as in \cite{ami},
then $Q_m(u_i)=0$ for each $m>i$. 
The set $\{ u_i: i \} $ is dense in $T$. Put 
$Y_1=C^0(T,{\bf K})$ and $Y_2=H$ 
and $J:=J_1\otimes J_2\in L_1(Y_1\otimes Y_2),$
where $J_1Q_m:=\alpha _mQ_m$ such that $\alpha _m\ne 0$ for each $m$
and $\sum_i|\alpha _i|<\infty $. Take $J_2$ also nondegenerate. Then $J$
induces a product measure $\mu $ on $C^0(T,H)$
such that $\mu =\mu _1\otimes \mu _2$, where $\mu _i$ are
measures on $Y_i$ induced by $J_i$ due to Formulas $(i,ii)$.
Analogously considering the following Banach subspace
$C^0_0(T,H):=\{ f\in C^0(T,H):$ $f(t_0)=0 \} $
and operators $J:=J_1\otimes J_2\in L_1(C^0_0(T,{\bf K})\otimes H)$
we get the measures $\mu $ on it also, where $t_0\in T$
is a marked point. 
\par For each finite number of points $(t_1,...,t_n)\subset T$ and
$(z_1,...,z_n)\subset H$ there exists a closed subset
$C^0(T,H;(t_1,...,t_n);(z_1,...,z_n)):=\{
f\in C^0(T,H):$ $f(t_i)=z_i;$ $i=1,...,n \} $ such that
$C^0(T,H;(t_1,...,t_n);(z_1,...,z_n))=(z_1,...,z_n)+
C^0(T,H;(t_1,...,t_n);(0,...,0))$, where
$C^0(T,H;(t_1,...,t_n);(0,...,0))$ is the Banach subspace 
of finite codimension $n$ in $C^0(T,H)$. Therefore, 
\par $(iii)$ $\sigma $-algebras $F_{t_2,t_1}^{-1}(Bf(H))$ and
$F_{t_4,t_3}^{-1}(Bf(H))$ are independent subalgebras
in the Borel $\sigma $-algebra $Bf(C^0_0(T,H))$, when $(t_1,t_2)$ and 
$(t_3,t_4)$ satisfy Condition $4.2.(i)$. 
\par Put $P(t_1,x_1,t_2,A):=\mu ( \{ f: f(t_1)=x_1, f(t_2)\in A \} )$
for each $t_1\ne t_2\in T,$ $x_1\in H$ and $A\in Bf(H)$. In view of $(iii)$
we get, that $P$ satisfies Conditions $3.2.(i-iv).$ By the above construction
(and Proposition 3.3.1 also) the Markov quasimeasure ${\tilde \mu }_{x_0}$
induced by $\mu $ is bounded, since $\mu $ is bounded, where
$x_0=0$ for $C^0_0(T,H)$. Let $\Omega $ be a set of elementary events
$\omega :=\{ f: f\in C^0_0(T,H), f(t_i)=x_i, i\in \Lambda _{\omega } \} ,$
where $\Lambda _{\omega }$ is a countable subset of $\bf N$,
$x_i\in H$, $(t_i: i\in \Lambda _{\omega } )$ is a subset of $T$
of pairwise distinct points. There exists the algebra $\tilde {\sf U}$
of cylindrical subsets of $C^0_0(T,H)$ induced by projections
$\pi _s: C^0_0(T,H)\to H^s,$ where $H^s:=\prod_{t\in s}H_t,$
$s=(t_1,...,t_n)$ are finite subsets of $T$, $H_t=H$ for each $t\in T$.
In view of the Kolmogorov theorem \cite{dalf,oksen,lu6,lubkt} 
$\tilde \mu _{x_0}$
on $((C^0_0(T,H),\tau _w),\tilde {\sf U})$ induces the probability measure
$\lambda $ on $(\Omega ,Bf(\Omega ))$, where $\tau _w$ is the weak topology
in $C^0_0(T,H)$.
\par Therefore, using product of measures we get examples of such measures $\mu $
for which stochastic processes exist
(see also Theorem 3.23, Lemmas 2.3, 2.5, 2.8 and \S 3.30 in \cite{lu6}).
Hence to each such measure
on $C^0_0(T,H)$ there corresponds the stochastic process.
Considering all operators
$J:=J_1\otimes J_2\in L_1(Y_1\otimes Y_2)$ and the corresponding 
measures as above we get ${\sf c}^{\aleph _0}=\sf c$ 
inequivalent measures by the Kakutani theorem II.4.1 \cite{dalf}
for each chosen $f$.
\par {\bf Note.} Evidently, this theorem is also true for $C^0(T,H)$, that
follows from the proof. If to take $\nu $ with $supp (\nu )=B({\bf K},0,1),$
then repeating the proof it is possible to construct $\mu $ with
$supp ( \mu )\subset B(C^0(T,{\bf K}),0,1)\times B(H,0,1).$
In the weak topology inherited from $C^0(T,H)$ the set $B(C^0(T,{\bf
K}),0,1)\times B(H,0,1)$ is compact and the condition $J\in L_1$ 
may be dropped. Certainly such measure $\mu $ 
can not be quasi-invariant relative
to shifts from a dense $\bf K$-linear subspace in $C^0(T,H)$,
but it can be constructed quasi-invariant relative to a dense
additive subgroup $G'$ of $B(C^0(T,{\bf K}),0,1)\times B(H,0,1)$, 
moreover, there exists $\mu $ for which $G'$ is also
$B({\bf K},0,1)$-absolutely convex.
\par {\bf 4.4.} We consider stochastic processes 
$E\in L^r(\Omega ,{\sf F},\lambda ;C^0(T,L_v(H)))$ 
such that $E=E(t,\omega ),$ where $1\le v\le \infty $, 
$1\le r\le \infty $, $t\in T=B({\bf K},t_0,R)$ and $\omega \in \Omega $ 
(see \S 2.14 and \S 4.2).
\par {\bf Definition.} For $L^r(\Omega ,{\sf F},\lambda ;C^0(T,L_v(H)))$ 
the non-Archimedean stochastic integral is defined by the following equation:
$$(i)\mbox{ }{\sf I}(E)(t,\omega ):=
(\hat P_wE)(t,\omega )=\sum_{j=0}^{\infty }E(t_j,\omega )[w(t_{j+1},\omega )
-w(t_j,\omega )],$$
where $w=w(t,\omega ),$ $t_j=\sigma _j(t)$ (see \S 2.1). 
\par {\bf 4.5.} {\bf Proposition.} {\it The non-Archimedean stochastic 
integral is the continuous $\bf K$-bilinear operator from 
$L^r(\Omega ,{\sf F},\lambda ;C^0(T,L_v(H)))
\otimes L^q(\Omega ,{\sf F},\lambda ;C^0_0(T,H))$ into
$L^s(\Omega ,{\sf F},\lambda ;C^0(T,H))$, where $1/q+1/r=1/s$
and $1\le r,q,s \le \infty $.} 
\par {\bf Proof.} It follows from Theorem  2.14,
since $(\hat P_{aw+by}E)=(a\hat P_wE)+(b\hat P_yE)$ 
and $(\hat P_w(aE+bV))=(a\hat P_wE)+b(\hat P_wV)$
for each $a, b\in \bf K$, each $w, y\in
L^q(\Omega ,{\sf F},\lambda ;C^0_0(T,H))$ and each
$E, V \in L^r(\Omega ,{\sf F},\lambda ;C^0(T,L_v(H)))$.
\par {\bf 4.6.} Consider a function 
$f$ from $T\times H$ into $Y=c_0(\beta ,{\bf K})$ satisfying conditions:
\par $(a)$ $f\in C^1(T\times H,Y);$
\par $(b)$ $(\bar \Phi ^nf)(t,x;h_1,...,h_n;\zeta _1,...,\zeta _n)
\in C^0(T\times H^{n+1}\times {\bf K^n},Y)$ for each $n\le m,$ 
\par $(c)$ $(\bar \Phi ^nf)(t,x;h_1,...,h_n;\zeta _1,...,\zeta _n)=0$ 
for $n=m+1$, 
\par $(d)$ $f(t,x)-f(0,x)=(\hat P_tg)(t,x)$ with $g\in C^0(T\times H,Y),$
where $2\le m\in \bf N$, $f=f(t,x)$, $t\in T$, $x\in H$; $h_1,...,h_n\in H$, 
$\zeta _1,...,\zeta _n\in \bf K;$ $\hat P_u$ is the antiderivation
operator on $C^0(T,Y)$, $(\hat P_tg)(t,x)$ is defined
for each fixed $x\in H$ by $t\in T$ such that
$(\hat P_tg)(t,x)=\hat P_ug(u,x)|_{u=t}$ with $u\in T$ 
(see \S 2.1 and also about difference quotients $(\bar \Phi ^nf)$
and spaces of functions of smoothness class $C^n$ in \cite{luseamb,lubp2}).
\par Suppose $a\in L^s(\Omega ,{\sf F},\lambda ;C^0(T,H)),$
$w\in L^q(\Omega ,{\sf F},\lambda ;C^0_0(T,H))$
and $E\in L^r(\Omega ,{\sf F},\lambda ;C^0(T,L(H)))$,
where $1/r+1/q=1/s$, $1\le r,q,s\le \infty $,
$a=a(t,\omega )$, $E=E(t,\omega )$, $t\in T,$ $\omega \in \Omega .$
A stochastic process of the type
$$(i)\quad \xi (t,\omega )=\xi _0(\omega )+
(\hat P_ua)(u,\omega )|_{u=t}+(\hat P_{w(u,\omega )}E)(
u,\omega )|_{u=t}$$
is said to have a stochastic differential
\par $(ii)$ $d\xi (t,\omega )=a(t,\omega )dt+E(t,\omega )dw(t,\omega )$,
since $(\hat P_tg)'(t)=g(t)$ for each $g\in C^0(T,H)$, where
$\xi _0\in L^s(\Omega ,{\sf F},\lambda ;H)$, $t_0, t\in T$,
$w(t_0,\omega )=0$. In view of Lemma 2.3, Theorem
2.14 and Proposition 4.5 $\xi \in L^s(\Omega ,{\sf F},\lambda ;C^0(T,H)).$
\par Let ${\hat P}_{u^b,w^h}$ denotes the antiderivation
operator ${\hat P}_{(\xi _1,...,\xi _{b+h})}$ given by Formula $2.1.(4)$,
where $\xi _1=u$,...,$\xi _b=u$, $\xi _{b+1}=w$,...,$\xi _{b+h}=w$.
Henceforth, it is used the notation
$$(iii) \mbox{ } {\tilde P}^n_{a,Ew}f(u,\xi (u,\omega )):=
\sum_{k=1}^n(k!)^{-1}\sum_{l=0}^k{k\choose l}(
{\hat P}_{u^{k-l},w(u,\omega )^l}$$
$$[(\partial ^kf/\partial x^k)
(u,\xi (u,\omega ))\circ (a^{\otimes (k-l)}\otimes E^{\otimes l})])$$
for such operator, when it exists (see the conditions above and below),
where $n\in \bf N$ or $n=\infty $.
\par {\bf Theorem.} {\it Let Conditions $4.6.(a-d),(i,ii)$ be satisfied, then
$$(iv)\mbox{ }f(t,\xi (t,\omega ))=f(t_0,\xi _0)+
(\hat P_uf'_t(u,\xi (u,\omega ))|_{u=t}+{\tilde P}^m_{a,Ew}f(u,\xi (u,\omega ))
|_{u=t}.$$}
\par {\bf Proof.} Let $\{ u_k:$ $k=0,1,...,n \} $ be a finite $|\pi |^l$
net in $T$, that is, for each $t\in T$ there exists
$k$ such that $|u_k-t|\le |\pi |^l$, where $n=n(k)\in \bf N$,
$\pi \in \bf K$, $p^{-1}\le |\pi |<1$ and $|\pi |$ is the generator
of the valuation group of $\bf K$,
since the ball $T$ is compact. We choose $t=u_n$ and $t_0=u_0$.
Denote by $\eta (t)$ a stochastic process $f(t,\xi (t,\omega ))$.
Then by the Taylor formula (see Theorem 29.4 \cite{sch1}
and Theorem 2.9 \cite{lupr180}) \\
$$(v)\quad f(t,\xi (t))-f(u,\xi (u))=
f'_t(u,\xi (u))(t-u)+f'_x(u,\xi (u)).(\Delta \xi )+$$
$$(1/2)f"_{t,t}(u,\xi (u))(t-u)^2
+f"_{t,x}(u,\xi (u)).((t-u),\Delta \xi )
+(1/2)f"_{x,x}(u,\xi (u)).(\Delta \xi ,$$ $$\Delta \xi )
+\{ (\bar \Phi ^2f)(u,\xi (u);(t-u),(t-u);1,1)
-(1/2)f"_{t,t}(u,\xi (u))(t-u)^2 \} $$
$$+\{ (\bar \Phi ^2f)(u,\xi (u);(t-u),\Delta \xi ;1,1)
+(\bar \Phi ^2f)(u,\xi (u);\Delta \xi ,(t-u);1,1)
-f"_{t,x}(u,\xi (u)$$ $$).(t-u,\Delta \xi ) \} 
+\{ (\bar \Phi ^2f)(u,\xi (u);\Delta \xi ,\Delta \xi ;1,1)
-(1/2)f"_{x,x}(u,\xi (u)).(\Delta \xi ,\Delta \xi ) \} ,$$
where $\Delta \xi =\xi (t)-\xi (u)$, for a brevity 
we denote $\xi (t)=\xi (t,\omega )$ and $w(t):=w(t,\omega )$
for a chosen $\omega $.
If $t_n=\sigma _n(t)$ for each $n=0,1,2,...$, then
by Formulas $(i)$ and 2.1.(4):
$$(vi)\mbox{ }\xi (t_{n+1},\omega )-\xi (t_n,\omega )=
a(t_n,\omega )(t_{n+1}-t_n)
+E(t_n,\omega )(w(t_{n+1},\omega )-w(t_n,\omega )),$$
where $ \{ \sigma _n: n=0,1,2,... \} $ is the approximation of the identity
in $T$.
\par From Condition $(d)$ it follows that 
$(\partial f(t,x)/\partial t)=g(t,x)={(\hat P_tg)'}_t$ and
$\hat P_t({f'}_t)(t,x)=f(t,x)-f(0,x)$,
which also leads to dissappearance
of terms $\partial ^{m+b}f(t,x)/\partial t^b\partial x^m$
from Formula $(iv)$ 
for each $b, m$ such that $1\le b$ and $2\le m+b$.
Now we approximate $f(t,x)$ by functions of the form
$\sum_j\phi _j(t)\psi _j(x)$, so the problem reduces
to the consideration of functions $f(x)$ which are independent
from $t$. Due to Conditions $(i,ii)$
it is possible to put $\xi (t,\omega )=
\xi _0(\omega )+a(\omega )(t-t_0)+E(\omega )[w(t)-w(t_0)]$. 
By the Taylor formula:
$$(vii) \mbox{ }f(x)=f(x_0)+\sum_{n=1}^m(n!)^{-1}f^{(
n)}(x_0).(x-x_0)^{\otimes n}$$
for each $x, x_0\in H,$ since $\bar \Phi ^{m+1}f=0$.
Put $t_k=\sigma _k(t)$ for each $k=0,1,2,...$, 
then $\eta (t)-\eta (t_0)=\sum_{j=0}^{\infty }
\{ f(\xi _{j+1}) -f(\xi _j) \} $, where $\xi _j:=\xi (t_j)$,
since $\lim_{j\to \infty }\xi _j=\xi $.
Then each term $f(\xi _{j+1}) -f(\xi _j)$ can be expressed by 
Formula $(vii)$ due to Condition $(b)$. 
On the other hand, $(\xi _{j+1}- \xi _j)=
a(\omega )(t_{j+1}-t_j)+E(\omega )[w(t_{j+1})-w(t_j)]$
as the particular case of Formula $(vi)$.
From Formulas $2.1.(4)$, $(v-vii)$ and Theorem 2.14 we get the statement 
of this theorem.
\par {\bf 4.7.} {\bf Corollary.} {\it If Conditions $4.6(a,d,i,ii)$
are satisfied, $4.6(b)$ is accomplished for each $n\in \bf N$ and
\par $(c')$ $\lim_{n\to \infty }
\|(\bar \Phi ^n_xf)(t,x;h_1,...,h_n;\zeta _1,...,\zeta _n)
\|_{C^0(T\times (B(H,0,R_1))^{n+1}
\times B({\bf K^{n+1}},0,R_1),Y)}=0$ for each $0<R_1<\infty ,$ then
$$(i)\mbox{ }f(t,\xi (t,\omega ))=f(t_0,\xi _0)+
(\hat P_uf'_t(u,\xi (u,\omega ))|_{u=t}+
({\tilde P}^{\infty }_{a,Ew}f(u,x))|_{u=t}.$$}
\par {\bf Proof.} From the proof of Theorem 4.6 we get a function $f(x)$
for which
$$(ii)\mbox{ }f(x)=f(x_0)+\sum_{n=1}^{\infty }(n!)^{-1}f^{(
n)}(x).(x-x_0)^{\otimes n}$$
due to Condition $(c')$. In view of Theorem 2.14
$$\lim_{m\to \infty }
\| (m!)^{-1}\sum_{l=0}^m{m\choose l}(
\hat P_{u^{m-l},w(u,\omega )^l}[(\partial ^mf/\partial x^m)
(u,\xi (u,\omega ))\circ $$
$$(a^{\otimes (m-l)}\otimes E^{\otimes l})])
|_{u=t}\| _{L^s(\Omega ,{\sf F},\lambda ;C^0(T,Y))}=0. $$
Approximating $f(x)$ by the Taylor formula up to terms ${\bar \Phi }^mf$
by finite sums and taking the limit while $m$ tends to the infinity
one deduces Formula $(i)$ from Formula $4.6(iv)$, since
for each chosen $\omega \in \Omega $ functions $a(t,\omega )$
and $w(t,\omega )$ are bounded on the compact ball $T$.
\par {\bf 4.8. Theorem.} {\it Let $f(u,x)\in C^{\infty }(T\times H,Y)$ 
and 
$$(i)\mbox{ } \lim_{n\to \infty }
\max_{0\le l\le n}\|(\bar \Phi ^nf)(t,x;h_1,...,h_n;$$
$$\zeta _1,...,\zeta _n)
\|_{C^0(T\times B({\bf K},0,r)^l\times B(H,0,1)^{n-l}
\times B({\bf K},0,R_1)^{n-l},Y)}=0$$
for each $0<R_1<\infty ,$ where $h_j=e_1$ and
$\zeta _j\in B({\bf K},0,r)$ for variables corresponding to 
$t\in T=B({\bf K},t_0,r)$
and $h_j\in B(H,0,1)$, $\zeta _j\in B({\bf K},0,R_1)$
for variables corresponding to $x\in H$, then
$$(ii)\quad f(t,\xi (t,\omega ))=f(t_0,\xi _0)+
\sum_{m+b\ge 1, 0\le m\in {\bf Z}, 0\le b\in {\bf Z}}
((m+b)!)^{-1}\sum_{l=0}^m{{m+b}\choose m} {m\choose l}$$
$$({\hat P}_{u^{b+m-l},w(u,\omega )^l}[(\partial ^{(m+b)}f/
\partial u^b\partial x^m)
(u,\xi (u,\omega ))\circ (I^{\otimes b}\otimes 
a^{\otimes (m-l)}\otimes E^{\otimes l})])
|_{u=t}.$$}
\par {\bf Proof.} In view of the Taylor formula we have
(see \cite{lupr180,sch1,sch2})
$$(iii)\quad f(t,x)=f(t_0,x_0)+\sum_{m+b=1}^k((m+b)!)^{-1}{{m+b}\choose m}
(\partial ^{(m+b)}f/\partial u^b\partial x^m)(t_0,x_0)$$
$$(t-t_0)^b.(x-x_0)^{\otimes m}+\sum_{m+b=k+1}{{k+1}\choose m}
[({\bar \Phi }^{k+1}f)(t_0,x_0;(t-t_0)^{\otimes b},
(x-x_0)^{\otimes m}; 1^{\otimes (k+1)})$$
$$-((k+1)!)^{-1}(\partial ^{(k+1)}f/\partial u^b\partial x^m)
(t_0,x_0)(t-t_0)^b.(x-x_0)^{\otimes m}]$$  for each $k\in \bf N$.
In view of Condition $(i)$, Formulas $(iii)$,
$2.1.(4)$, $4.6.(vi)$ we get Formula $(ii)$ (see the proof of Theorem 4.6).


\begin{thebibliography}{99}
\bibitem{aiel} S. Aida, D. Elworthy. "Differential calculus
on path and loop spaces. 1. Logarithmic Sobolev inequalities
on path spaces". C.R. Acad. Sci. Paris. Ser. 1, {\bf 321} (1995), 
97-102.
\bibitem{ami} Y. Amice. "Interpolation p-adique". Bull. Soc. Math. France
{\bf 92}(1964), 117-180.
\bibitem{bikvol} A.H. Bikulov, I.V. Volovich.
"$p$-adic Brownian motion". Izv. Ross. Akad. Nauk.
Ser. Math. {\bf 61: 3} (1997), 75-90.
\bibitem{boui}
N. Bourbaki. "Integration". Chapters 1-9 (Moscow: Nauka, 1970 and 1977).
\bibitem{dalf} Yu.L. Dalecky, S.V. Fomin. "Measures and differential
equations in infinite-dimensional space" (Dordrecht, The Netherlands:
Kluwer, 1991).
\bibitem{dalschn} Yu.L. Dalecky, Ya.I. Schnaiderman.
"Diffusion and quasi-invariant measures on infinite-dimensional
Lie groups". Funct. Anal. and Pril. {\bf 3: 2} (1969), 88-90.
\bibitem{eljm} K.D. Elworthy, Y. Le Jan, X.-M. Li.
"Integration by parts formulae for degenerate diffusion measures
on path spaces and diffeomorphism groups". C.R. Acad. Sci. Paris.
Ser. 1, {\bf 323} (1996), 921-926.
\bibitem{eng} R. Engelking. "General topology" (Moscow: Mir,1986).
\bibitem{esc} A. Escassut. "Analytic elements in $p$-adic analysis"
(Singapore: World Scientific, 1995).
\bibitem{evans}  S.N. Evans. "Continuity properties
of Gaussian stochastic processes indexed by a local field".
Proceed. Lond. Math. Soc. Ser. 3, {\bf 56} (1988), 380-416.
\bibitem{feder} H. Federer. "Geometric measure theory"
(Berlin: Springer, 1968).
\bibitem{fell} J.M.G.  Fell, R.S.  Doran.  "Representations
of $*$-Algebras, Locally Compact Groups, and Banach $*$-Algebraic Bundles"
V. {\bf 1} and V. {\bf 2}.  (Acad.  Press, Boston, 1988).  
\bibitem{freput} J. Fresnel, M. van der Put. "G\'eom\'etrie
analytique rigide et applications" (Boston, Birkh\"auser, 1981).
\bibitem{gihsko} I.I. Gihman, A.V. Skorohod. "Stochastic differential
equations and their apllications" (Kiev: Naukova Dumka, 1982).
\bibitem{grus} L. Gruson. "Th\'eorie de Fredholm $p$-adique".
Bull. Soc. Math. France. {\bf 94} (1966), 67-95.
\bibitem{hew} E. Hewitt, K.A. Ross. "Abstract harmonic analysis"
(Berlin: Springer, 1979).
\bibitem{ikwat} N. Ikeda, S. Watanabe. "Stochastic differential equations
and diffusion processes" (Moscow: Nauka, 1986).
\bibitem{jang} Y. Jang. "Non-Archimedean quantum mechanics".
Tohoku Mathem. Publications. $N^o$ {\bf 10} 
(Tohoku: Toh. Univ., Math. Inst., 1998).
\bibitem{khrif} A.Yu. Khrennikov. "Generalized functions and Gaussian 
path integrals". Izv. Acad. Nauk. Ser. Mat. {\bf 55} 
(1991), 780-814.
\bibitem{lud} S.V. Ludkovsky. "Measures on groups of diffeomorphisms of
non-Archimedean Banach manifolds". 
Russ. Math. Surveys. {\bf 51: 2} (1996), 338-340.
\bibitem{lupr180} S.V. Ludkovsky. "Representations and structure
of groups of diffeomorphisms of non-Archimedean Banach manifolds".
Intern.  Centre for Theoret.  
Phys., Trieste, Italy, Preprint {\bf IC/96/180}, 
23 pages, September, 1996.
\bibitem{lu6} S.V. Ludkovsky. "Quasi-invariant 
and pseudo-differentiable measures on a non-Archimedean Banach space";
Intern.  Centre for Theoret.  
Phys., Trieste, Italy, Preprints
{\bf IC/96/210}, 50 pages, October, 1996.  
\bibitem{lubkt} S.V. Ludkovsky. 
"The non-Archimedean analogs of the Bochner-Kolmogorov,
Minlos-Sazonov and Kakutani theorems".
Los Alamos National Laboratory, USA.
Preprint {\bf math.FA/0010230}, 32 pages, 25 October 2000.
\bibitem{luumnls} S.V. Ludkovsky. "Quasi-invariant measures on
non-Archimedean semigroups of loops". Russ. Math. Surveys {\bf 53: 3} (1998),
633-634.
\bibitem{luseamb} S.V. Ludkovsky. "Irreducible unitary representations
of non-Archimedean groups of diffeomorphisms". Southeast Asian
Mathem. Bull. {\bf 22: 3} (1998), 301-319. 
\bibitem{lubp99} S.V. Ludkovsky. "Properties of quasi-invariant measures on
topological groups and associated algebras".
Annales Math\'ematiques Blaise Pascal. {\bf 6: 1} (1999), 33-45.
\bibitem{lutmf99} S.V. Ludkovsky. "Measures on groups of diffeomorphisms
of non-Archimedean manifolds, representations of groups and their 
applications". Theoret. and Math. Phys. {\bf 119: 3} (1999), 698-711.
\bibitem{lubp2} S.V. Ludkovsky. "Quasi-invariant measures on non-Archimedean 
groups and semigroups of loops and paths, their representations. I, II".
Annales Math\'ematiques Blaise Pascal. {\bf 7: 2} (2000), 19-53, 55-80.
\bibitem{lustrcm} S.V. Ludkovsky.
"Stochastic processes on geometric loop groups and diffeomorphism 
groups of real and complex manifolds, associated unitary representations",
Los Alamos National Laboratory, USA.
Preprint {\bf math.GR/0102222}, 35 pages, 28 February 2001.
\bibitem{ludgla} S.V. Ludkovsky.
"A structure and representations of diffeomorphism groups
of non-Archimedean manifolds".
Infinite dimensional analysis, quantum probability and related topics. 
To appear (previous variant: Los Alamos National Laboratory, USA.
Preprint {\bf math.GR/0004126}, 32 pages, 19 April 2000).
\bibitem{mckean} H.P. Mc Kean. "Stochastic integrals" 
(Moscow: Mir, 1972).
\bibitem{malli} P. Malliavin. "Stochastic analysis" 
(Berlin: Springer, 1997).
\bibitem{nari} L. Narici, E. Beckenstein. "Topological vector spaces"
(New York: Marcel Dekker Inc., 1985).
\bibitem{oksen} B. {\O}ksendal. "Stochastic differential equations"
(Berlin: Springer, 1995).
\bibitem{pietsch} A. Pietsch. "Nukleare lokalkonvexe R\"aume"
(Berlin: Akademie-Verlag, 1965).
\bibitem{put} M. van der Put. "The ring of bounded operators
on a non-Archimedean normed linear space".
Indag. Math. {\bf 71: 3} (1968), 260-264.
\bibitem{reed} M. Reed, B. Simon. "Methods of modern mathematical physics"
(New York: Acad. Press, 1975).
\bibitem{roo} A.C.M. van Rooij. "Non-Archimedean functional analysis"
(New York: Marcel Dekker Inc., 1978).
\bibitem{sch1} W.H. Schikhof. "Ultrametric calculus" (Cambridge:
Cambr. Univ. Press, 1984).  
\bibitem{sch2} W.H. Schikhof. "Non-Archimedean calculus". Nijmegen: Math.
Inst., Cath. Univ., Report {\bf 7812}, 130 pages, 1978.
\bibitem{schikop} W.H. Schikhof. "On $p$-adic compact operators".
Report 8911 (Dep. Math. Cath. Univ., Nijmegen, The Netherlands, 1989).
\bibitem{vla3} V.S. Vladimirov, I.V. Volovich, E.I. Zelenov. "$p$-Adic
analysis and mathematical physics" (Moscow: Fiz.-Mat. Lit, 1994).
\bibitem{wei} A. Weil. "Basic number theory" (Berlin: Springer, 1973).
\end{thebibliography}
\end{document}